\documentclass[brochure,12pt]{bourbaki}
\usepackage[matrix,arrow]{xy}
\usepackage{amssymb,amsfonts,amsmath,footnote}
\usepackage{enumerate}
\usepackage[francais]{babel}
\usepackage{comment}
\usepackage{footmisc}
\usepackage[latin1]{inputenc}
\DeclareFontFamily{OT1}{pzc}{}
\DeclareFontShape{OT1}{pzc}{m}{it}{<-> s * [1.10] pzcmi7t}{}
\DeclareMathAlphabet{\mathpzc}{OT1}{pzc}{m}{it}

\newtheorem{ques}[defi]{\quesname}
 \def\quesname{Question}

\DeclareMathOperator{\GL}{GL}

\DeclareMathOperator{\Frac}{Frac}
\DeclareMathOperator{\LEnd}{LEnd}
\DeclareMathOperator{\Spec}{Spec}
\DeclareMathOperator{\dR}{dR}
\DeclareMathOperator{\B}{B}
\DeclareMathOperator{\Id}{Id}
\DeclareMathOperator{\End}{End}
\DeclareMathOperator{\lisse}{lisse}
\DeclareMathOperator{\Ab}{Ab}
\DeclareMathOperator{\GSp}{GSp}
\DeclareMathOperator{\Cl}{Cl}
\DeclareMathOperator{\Sh}{Sh}
\DeclareMathOperator{\SO}{SO}
\DeclareMathOperator{\Spin}{Spin}
\DeclareMathOperator{\GSpin}{GSpin}

\DeclareMathOperator{\prim}{prim}
\DeclareMathOperator{\cris}{cris}
\DeclareMathOperator{\et}{\mathrm{\acute{e}t}}
\DeclareMathOperator{\Gal}{Gal}
\DeclareMathOperator{\Pic}{Pic}
\DeclareMathOperator{\Ker}{Ker}
\DeclareMathOperator{\Br}{Br}
\DeclareMathOperator{\Spf}{Spf}
\DeclareMathOperator{\Mor}{Mor}
\DeclareMathOperator{\ovQ}{\overline{\mathbb{Q}}}
\DeclareMathOperator{\ovFp}{\overline{\mathbb{F}}_{\mathit{p}}}

\date{Mars 2014}
\bbkannee{66\`eme ann\'ee, 2013-2014}
\bbknumero{1081}
\title{CONSTRUCTION DE COURBES SUR LES SURFACES K3}
\subtitle{d'apr\`es Bogomolov-Hassett-Tschinkel, Charles,  Li-Liedtke, \\Madapusi Pera, Maulik...}
\author{Olivier BENOIST}
\address{CNRS, IRMA\\
Universit\'e de Strasbourg\\
7, rue Ren\'e Descartes\\
F--67000 Strasbourg}
\email{obenoist@unistra.fr}

\begin{document}
\maketitle

\noindent{\bf INTRODUCTION}

On expose ici des r\'esultats r\'ecents de construction de courbes sur les surfaces K3 : la preuve de la conjecture de Tate pour les surfaces K3 en caract\'eristique finie impaire (Th\'eor\`eme \ref{thTate2}, d'apr\`es Maulik \cite{Maulik}, Charles \cite{Charles} et Madapusi Pera \cite{MP}), et la construction d'une infinit\'e de courbes rationnelles sur de nombreuses surfaces K3 (Th\'eor\`eme \ref{BHTLL2}, d'apr\`es Bogomolov-Hassett-Tschinkel \cite{BHT} et Li-Liedtke \cite{LL}).

La premi\`ere partie de ce texte est consacr\'ee \`a l'\'etude de l'application de Kuga-Satake, n\'ecessaire \`a la preuve de la conjecture de Tate dans la deuxi\`eme partie. Finalement, la troisi\`eme partie exploite la conjecture de Tate pour construire des courbes rationnelles sur des surfaces K3. Les introductions de chacune de ces parties d\'ecrivent plus pr\'ecis\'ement leur contenu.

Nous rappelons ici des g\'en\'eralit\'es sur les surfaces K3 en caract\'eristique finie, puis nous pr\'esentons les r\'esultats mentionn\'es ci-dessus.

\begin{defi}
Une \textbf{surface K3} sur un corps $k$ est une surface projective\footnote{En particulier, quand $k=\mathbb{C}$, on ne consid\'erera pas dans ce texte de surfaces K3 non alg\'ebriques.} lisse et g\'eom\'etriquement connexe $X$ sur $k$ telle que $K_X\simeq\mathcal{O}_X$ et $H^1(X,\mathcal{O}_X)=0$.

Les classes d'isomorphisme de fibr\'es en droites sur $X$ forment un groupe ab\'elien libre de type fini $\Pic(X)$. On note $\rho(X)$ son rang : c'est le \textbf{nombre de Picard} de $X$.
\end{defi}

\subsection{Surfaces K3 complexes}

Si $X$ est une surface K3 complexe, l'injectivit\'e de l'application classe de cycle en cohomologie de Betti $\Pic(X)\to H^2(X,\mathbb{Z})$ montre que $\rho(X)\leq b_2(X)=22$. La th\'eorie de Hodge implique que l'image de cette application est  incluse dans $H^{1,1}(X)$, ce qui fournit l'in\'egalit\'e plus forte $\rho(X)\leq 20$. Enfin, le th\'eor\`eme des classes $(1,1)$ de Lefschetz montre que $\Pic(X)\to H^2(X,\mathbb{Z})\cap H^{1,1}(X)$ est un isomorphisme, et permet de calculer $\rho(X)$ connaissant la structure de Hodge de $X$.

Le nombre de Picard peut prendre toutes les valeurs entre $1$ et $20$, et la th\'eorie des d\'eformations montre que, dans un espace des modules des surfaces K3 polaris\'ees, le lieu o\`u $\rho(X)\geq r$ est une r\'eunion d\'enombrable de sous-vari\'et\'es de dimension $20-r$.

\subsection{Surfaces K3 de hauteur finie}

  Soit maintenant $X$ une surface K3 sur un corps $k$ alg\'ebriquement clos de carac\-t\'eristique finie. Igusa \cite{Igusa} a montr\'e que l'in\'egalit\'e 
$\rho(X)\leq b_2(X)=22$ \'etait encore v\'erifi\'ee. Pour obtenir des obstructions suppl\'ementaires \`a l'existence de fibr\'es en droites, analogues \`a celles obtenues sur $\mathbb{C}$ par th\'eorie de Hodge, Artin et Mazur \cite{ArtinMazur} ont introduit un nouvel invariant des surfaces K3. Ils consid\`erent le foncteur contravariant $T\mapsto\Ker[ \Br(X_T)\to \Br(X)]$ d\'efini sur les $k$-sch\'emas finis locaux. Ils montrent que ce foncteur est repr\'esentable par un groupe formel lisse de dimension $1$ sur $k$ : le \textbf{groupe de Brauer formel} $\widehat{\Br}(X)$. Comme ces groupes formels sont classifi\'es  par leur hauteur \cite{Lazard}, on obtient un invariant $h(X)\in \mathbb{N}^*\cup\{\infty\}$ : c'est la \textbf{hauteur} de $X$. Cet invariant varie semi-contin\^ument sup\'erieurement avec $X$.

  Si $X$ est une surface K3 de hauteur finie, Artin et Mazur \cite{ArtinMazur} montrent que \mbox{$\rho(X)\leq 22-2h(X)$.} Ces surfaces K3 se comportent comme les surfaces K3 en carac\-t\'eristique nulle : elles v\'erifient $\rho(X)\leq 20$ et, si l'on se restreint aux surfaces K3 de hauteur finie, le lieu dans un espace de modules de surfaces K3 polaris\'ees o\`u $\rho(X)\geq r$ est une r\'eunion d\'enombrable de sous-vari\'et\'es de dimension $20-r$.

  
Une surface K3 de hauteur $1$ est dite \textbf{ordinaire}, et les surfaces K3 ordinaires forment un ouvert dense des espaces de modules de surfaces K3 polaris\'ees sur $k$.

\subsection{Surfaces K3 supersinguli\`eres} 

Les surfaces K3 de hauteur infinie (i.e. telles que $\widehat{\Br}(X)\simeq\widehat{\mathbb{G}_a}$) sont dites \textbf{supersinguli\`eres} (ou \textbf{Artin-supersinguli\`eres}). Leur d\'efinition, de nature cohomologique, admet plusieurs formulations \'equivalentes :

\begin{prop}
\label{defss}
Les assertions  (i) et  (ii) sont \'equivalentes, et sont \'equivalentes \`a  (iii) si $k=\ovFp$
et $\ell$ est un nombre premier diff\'erent de $p$.
\begin{enumerate}[(i)]
\item $X$ est supersinguli\`ere.
\item Les pentes de $H^2_{\cris}(X/W(k))_{\mathbb{Q}}(1)$ sont toutes \'egales \`a $0$.
\item  L'endomorphisme de Frobenius de $H^2_{\et}(X,\mathbb{Q}_\ell(1))$ est d'ordre fini. 
\end{enumerate}
\end{prop}

\noindent{\sc Preuve} ---   L'\'equivalence des deux premiers \'enonc\'es est due \`a Artin et Mazur \cite{ArtinMazur}.
Si $X$ est d\'efinie sur un corps fini $\mathbb{F}_q$, Katz et Messing \cite{KM} montrent que les valeurs propres du Frobenius g\'eom\'etrique agissant en cohomologie $\ell$-adique et cristalline co\"incident.
Ces valeurs propres sont des unit\'es $\ell$-adiques par dualit\'e de Poincar\'e, leurs valuations \mbox{$p$-adiques} correctement normalis\'ees sont les pentes de $H^2_{\cris}(X/W(k))_{\mathbb{Q}}(1)$, et leurs valeurs absolues complexes sont \'egales \`a $1$ par les conjectures de Weil prouv\'ees par Deligne \cite{WeilK3}. L'assertion (ii) est donc \'equivalente au fait que ce sont des entiers alg\'ebriques dont tous les conjugu\'es complexes sont de module $1$. Ceci \'equivaut au fait que ce sont des racines de l'unit\'e, par un lemme de Kronecker, et \`a l'assertion (iii)
car l'action du Frobenius sur $H^2_{\et}(X,\mathbb{Q}_\ell(1))$ est semi-simple \cite[Corollaire 1.10]{K3rel}.\qed

\bigskip

Les propri\'et\'es des surfaces K3 supersinguli\`eres sont tr\`es diff\'erentes de celles des surfaces K3 complexes. On ne dispose pas d'obstruction suppl\'ementaire \`a l'existence de fibr\'es en droites : Tate \cite{Tatepoles} et Shioda \cite{Shiodass} ont donn\'e des exemples de telles surfaces avec nombre de Picard $22$ (dites \textbf{Shioda-supersinguli\`eres}). Artin a conjectur\'e que toutes les surfaces K3 supersinguli\`eres \'etaient Shioda-supersinguli\`eres, et on verra ci-dessous (corollaire \ref{cororho} (ii)) que cette conjecture est vraie au moins en caract\'eristique impaire.

Les surfaces K3 supersinguli\`eres constituent un ferm\'e de dimen\-sion $9$
des espaces de modules de surfaces K3 polaris\'ees sur $k$ \cite[Theorem 15]{O3}.

\subsection{La conjecture de Tate}

La conjecture de Tate, qu'on peut \'enoncer sur tout corps de type fini\footnote{On dit qu'un corps est de type fini s'il est engendr\'e par un nombre fini d'\'el\'ements sur son sous-corps premier.}, rend 
les m\^emes services que le th\'eor\`eme des classes $(1,1)$ de Lefschetz, en calculant $\rho(X)$ \`a l'aide de donn\'ees cohomologiques :

\begin{conj}
\label{Tate}
Soit $X$ une surface K3 sur un corps $k$ de type fini. Si $\ell$ est un nombre premier inversible dans $k$, et si $\bar{k}$ est une cl\^oture s\'eparable de $k$, l'application classe de cycle induit un isomorphisme :
$$\Pic(X)\otimes \mathbb{Q}_\ell\to H^2_{\et}(X_{\bar{k}},\mathbb{Q}_\ell(1))^{\Gal(\bar{k}/k)}.$$
\end{conj}

Faisant suite aux travaux classiques \cite{ASD, Nygaard, NO}, des progr\`es r\'ecents de Maulik \cite{Maulik}, Charles \cite{Charles} et Madapusi Pera \cite{MP} ont permis d'obtenir :

\begin{theo} 
\label{thTate1}
La conjecture \ref{Tate} est vraie en caract\'eristique diff\'erente de $2$.
\end{theo}

On renvoie \`a l'introduction de la deuxi\`eme partie de ce texte pour une discussion des contributions respectives de ces articles.

La conjecture \ref{Tate} est un cas particulier d'une conjecture bien plus g\'en\'erale, due \`a Tate \cite{Tatepoles}, qui pr\'edit l'existence de cycles alg\'ebriques en toute codimension sur toute vari\'et\'e projective lisse sur un corps de type fini. On invite le lecteur \`a consulter l'article de survol \cite{Srinivas}, qui replace notamment la conjecture \ref{Tate} dans ce cadre g\'en\'eral.

Tr\`es peu d'instances de cette conjecture sont connues. Si le cas des diviseurs sur les surfaces K3 est accessible, c'est gr\^ace \`a leur lien avec les vari\'et\'es ab\'eliennes, fourni par la construction de Kuga-Satake.
Pour cette raison, la premi\`ere partie de ce texte est consacr\'ee \`a l'\'etude de l'application de Kuga-Satake, et la deuxi\`eme partie \`a la preuve proprement dite du th\'eor\`eme \ref{thTate1}.

\subsection{Cons\'equences de la conjecture de Tate}

La conjecture de Tate permet d'obtenir des informations sur le nombre de Picard $\rho(X)$. Dans le corollaire ci-dessous, (i) \'etait d\'ej\`a connu d'Artin \cite{Artin} en toute carac\-t\'eristique, (ii) est la conjecture d'Artin mentionn\'ee ci-dessus et l'argument de (iii), attribu\'e \`a Swinnerton-Dyer, est \cite[Theorem 13]{BHT}.

\begin{coro}\label{cororho}
Soit $X$ une surface K3 sur un corps $k$ alg\'ebriquement clos de caract\'eristique $p\geq 3$. Alors :
\begin{enumerate}[(i)]
\item $\rho(X)\in \{1,\dots, 20,22\}$.
\item  $\rho(X)=22$ si et seulement si $X$ est supersinguli\`ere.
\item Si $X$ est d\'efinie sur un corps fini, $\rho(X)$ est pair.
\end{enumerate}
\end{coro}

\noindent{\sc Preuve} ---  Si $X$ n'est pas supersinguli\`ere, $\rho(X)\leq 22-2h(X)\leq 20$. Si $X$ est supersinguli\`ere, choisissons-en une sp\'ecialisation $Y$ d\'efinie sur un corps fini : $Y$ est encore supersinguli\`ere. Par la proposition \ref{defss} (iii), on peut supposer que l'action de Galois sur $H^2_{\et}(Y_{\ovFp},\mathbb{Q}_\ell(1))$ est triviale, et le th\'eor\`eme \ref{thTate1} montre que $\rho(Y_{\ovFp})=22$. Par un th\'eor\`eme d'Artin \cite[Theorem 1.1]{Artin}, $\rho(X)=\rho(Y_{\ovFp})=22$. Nous avons montr\'e (i) et (ii).

Supposons que $X$ est d\'efinie sur un corps fini. Quitte \`a \'etendre les scalaires, on peut supposer que la seule valeur propre du Frobenius agissant sur $H^2_{\et}(X_{\ovFp},\mathbb{Q}_\ell(1))$ qui soit une racine de l'unit\'e est $1$. Par dualit\'e de Poincar\'e, si $\alpha$ est une autre valeur propre, $\alpha^{-1}$ en est \'egalement une, avec m\^eme multiplicit\'e. On en d\'eduit que la multiplicit\'e de la valeur propre $1$ est paire. Comme l'action du Frobenius sur $H^2_{\et}(X_{\ovFp},\mathbb{Q}_\ell(1))$ est semi-simple \cite[Corollaire 1.10]{K3rel}, le th\'eor\`eme \ref{thTate1} montre (iii). \qed

\bigskip

De nombreux \'enonc\'es sur les surfaces K3 supersinguli\`eres n'\'etaient auparavant connus que pour les surfaces K3 Shioda-supersinguli\`eres. On peut mentionner les travaux d'Ogus sur le th\'eor\`eme de Torelli cristallin \cite{O1, O2}, et le th\'eor\`eme suivant de Rudakov et Shafarevich \cite{RS} :

\begin{theo}\label{RudaSha}
Une surface K3 supersinguli\`ere sur le corps des fractions d'un anneau de valuation discr\`ete de caract\'eristique $\geq 5$ a potentiellement bonne r\'eduction. 
\end{theo}

Enfin, Lieblich, Maulik et Snowden \cite{LMS} ont r\'ecemment montr\'e la cons\'equence suivante de la conjecture de Tate, analogue pour les surfaces K3 d'un th\'eor\`eme de Zarhin \cite{Zarhin} pour les vari\'et\'es ab\'eliennes :

\begin{theo}
Si $k$ est fini de caract\'eristique $\geq 5$, il n'existe qu'un nombre fini de classes d'isomorphisme de surfaces K3 sur $k$.
\end{theo}

\subsection{Construction de courbes rationnelles sur les surfaces K3}

Nous avons jusqu'ici discut\'e l'existence de classes de courbes sur les surfaces K3. On peut aussi chercher \`a d\'emontrer l'existence de courbes particuli\`eres. Par exemple, une \textbf{courbe rationnelle} est une sous-vari\'et\'e int\`egre dont la normalisation est isomorphe \`a $\mathbb{P}^1$. On conjecture :

\begin{conj}\label{conjcourbes}
Une surface K3 sur un corps alg\'ebriquement clos $k$ contient une infinit\'e de courbes rationnelles. 
\end{conj}

Insistons sur le fait que cette conjecture est d\'ej\`a int\'eressante, et toujours ouverte en g\'en\'eral, quand $k=\mathbb{C}$. En utilisant de mani\`ere cruciale la conjecture de Tate pour les surfaces K3 sur les corps finis, et plus pr\'ecis\'ement le corollaire \ref{cororho} (iii), Bogomolov, Hassett et Tschinkel \cite{BHT} et Li et Liedtke \cite{LL} y ont r\'epondu positivement dans de nombreux cas :

\begin{theo} 
\label{BHTLL}
La conjecture \ref{conjcourbes} est vraie dans les deux cas suivants :
\begin{enumerate}[(i)]
\item Si $k$ est de caract\'eristique nulle et $X$ n'est pas d\'efinie sur $\ovQ$.
\item Si $k$ est de caract\'eristique $\neq 2$ et $\rho(X)$ est impair.
\end{enumerate}
\end{theo}

La troisi\`eme partie de ce texte, ind\'ependante des deux premi\`eres, est consacr\'ee \`a la d\'emonstration de ce th\'eor\`eme ; on renvoie \`a son introduction pour plus d'informations.

\bigskip
{\bf Remerciements.}
Merci \`a Fran\c{c}ois Charles, \`a qui ce texte doit beaucoup, pour les innombrables discussions que nous avons eues sur les math\'ematiques abord\'ees ici ! Je suis aussi reconnaissant \`a
Keerthi Madapusi Pera, qui a r\'epondu avec pr\'ecision \`a mes nombreuses questions sur ses travaux, et \`a Davesh Maulik pour le plaisir que j'ai eu \`a lire \cite{Maulik}.
Je tiens enfin
\`a remercier H\'el\`ene Esnault, Javier Fres\'an, Daniel Huybrechts, Christian Liedtke et Gianluca Pacienza, qui ont relu des versions pr\'eliminaires de ce texte : leurs conseils ont beaucoup am\'elior\'e celui-ci.

\section{L'application de Kuga-Satake}

     Kuga et Satake \cite{KS} ont associ\'e \`a toute surface K3 complexe polaris\'ee $X$ une vari\'et\'e ab\'elienne $A$ : sa vari\'et\'e de Kuga-Satake. La construction, purement transcendente, consiste \`a munir l'alg\`ebre de Clifford paire de $H^2_{\prim}(X,\mathbb{Z})$ d'une structure de Hodge polaris\'ee de poids $1$, qui correspond \`a une vari\'et\'e ab\'elienne $A$.
   
     Des travaux successifs de Deligne \cite{WeilK3}, Ogus \cite{OgusII}, Andr\'e \cite{Andre} ont permis d'associer \`a une surface K3 sur tout corps  une vari\'et\'e de Kuga-Satake, et d'\'etendre \`a d'autres th\'eories cohomologiques ($\ell$-adique, cristalline) les liens qui sont visibles, dans la d\'efinition transcendente, entre les cohomologies de Betti de $X$ et de $A$. Ces techniques ont permis de nombreux progr\`es dans l'\'etude des surfaces K3 : le premier exemple en a \'et\'e la preuve par Deligne des conjectures de Weil pour les surfaces K3 \cite{WeilK3}.
   
   Les avanc\'ees sur la conjecture de Tate pour les surfaces K3, expliqu\'ees dans la partie suivante, s'inscrivent dans cette lign\'ee. Elles utilisent de mani\`ere cruciale un r\'esultat nouveau sur la construction de Kuga-Satake. On peut voir cette construction comme un avatar de l'application des p\'eriodes pour les surfaces K3 pour laquelle on dispose, en caract\'eristique nulle, d'un th\'eor\`eme de Torelli. Rizov \cite{RizovKS}, Maulik \cite{Maulik} et Madapusi Pera \cite{MP} ont \'etendu, avec une g\'en\'eralit\'e croissante, ce th\'eor\`eme de Torelli en caract\'eristique mixte. C'est l'objectif principal de cette partie (th\'eor\`eme \ref{Torelli}).
   
   L'\'etude de l'application de Kuga-Satake en caract\'eristique mixte a \'et\'e initi\'ee par Rizov \cite{RizovKS}. Nous suivons ici de tr\`es pr\`es la pr\'esentation de Madapusi Pera \cite{MP}, qui permet de travailler sur $\mathbb{Z}[\frac{1}{2}]$, et nous renvoyons \`a cet article pour plus de d\'etails.
   Comme Charles \cite{Charles} et Madapusi Pera \cite{MP}, nous utilisons une variante de la construction de Kuga-Satake classique, prenant en compte l'alg\`ebre de Clifford toute enti\`ere, qui simplifie les relations cohomologiques entre $X$ et $A$.

\subsection{La construction de Kuga-Satake}\label{KSnaive}

Si $X$ est une surface K3 polaris\'ee, on note $P^2(X):=H^2_{\prim}(X)$ la cohomologie primitive de $X$, et on la munit de la forme quadratique $\langle\cdot,\cdot\rangle$ induite par la dualit\'e de Poincar\'e \footnote{De nombreux auteurs (\cite{Andre}, \cite{RizovKS}, \cite{Maulik}, \cite{MP}) utilisent l'oppos\'e de cette forme quadratique.}. On utilisera des indices $\B,\dR, \et, \ell, \cris$ pour indiquer la th\'eorie cohomologique (Betti, de Rham, \'etale, $\ell$-adique, cristalline) utilis\'ee.

Soit $X$ une surface K3 complexe polaris\'ee, et soit $\omega=x+iy$ un g\'en\'erateur de $H^{2,0}(X)$ tel que $\langle \omega,\bar{\omega}\rangle=2$. 
Munissons l'alg\`ebre de Clifford r\'eelle $\Cl(P^2_{\B}(X,\mathbb{R}))$ d'une structure complexe en faisant agir $x\cdot y$ par multiplication \`a gauche.
Alors $A:=\Cl(P^2_{\B}(X,\mathbb{R}))/\Cl(P^2_{\B}(X,\mathbb{Z}))$ est un tore complexe. On montre (\cite{KS}, \cite[Chapter 4]{HuyK3}) qu'il s'agit d'une vari\'et\'e ab\'elienne : la \textbf{vari\'et\'e de Kuga-Satake} de $X$.

Les vari\'et\'es $X$ et $A$ sont li\'ees par un morphisme de structures de Hodge induit par la multiplication \`a droite de $P^2_{\B}(X,\mathbb{Z})$ sur $\Cl(P^2_{\B}(X,\mathbb{Z}))$ :\begin{equation}\label{KSrel}P^2_{\B}(X)(1)\hookrightarrow H^{\otimes(1,1)}_{\B}(A):=\End(H^1_{\B}(A)).\end{equation}

On note $L_{\B}(A)\subset H^{\otimes(1,1)}_{\B}(A)$ la sous-structure de Hodge image.

La relation (\ref{KSrel}) sera cruciale dans la suite : elle permettra d'\'etablir un lien entre fibr\'es en droites sur $X$ et endomorphismes de $A$ (voir par exemple la proposition \ref{ident}).

\subsection{Espaces de modules de surfaces K3}\label{edm}

Dans les paragraphes qui suivent, on explique comment effectuer cette construction en famille, sur d'autres bases que $\mathbb{C}$, et comment \'etendre la relation (\ref{KSrel}) \`a d'autres th\'eories cohomologiques.
  Commen\c{c}ons par introduire les espaces de modules dont nous aurons besoin.

Fixons un entier $d\geq 1$, et notons
$\mathfrak{M}_{2d}$ le champ de modules sur $\mathbb{Z}[\frac{1}{2}]$ des surfaces K3 munies d'une polarisation
primitive de degr\'e $2d$ \cite{RizovK3}.
C'est un champ de Deligne-Mumford s\'epar\'e de type fini sur $\mathbb{Z}[\frac{1}{2}]$.
On consid\'erera plut\^ot $\widetilde{\mathfrak{M}}_{2d}$
son rev\^etement double \'etale correspondant aux trivialisations isom\'etriques du d\'eterminant de la cohomologie primitive $2$-adique\footnote{Par \cite[Lemma 3.2]{Saitodisc}, le nombre $2$ ne joue pas ici un r\^ole particulier.}
des surfaces K3.

Il r\'esulte de \cite[Proposition 2.2]{O1} que le lieu non lisse de $\widetilde{\mathfrak{M}}_{2d}$ est constitu\'e de points isol\'es correspondant \`a des surfaces K3 supersinguli\`eres particuli\`eres, dites super\-sp\'eciales. De tels points n'existent qu'en caract\'eristique $p$ avec $p\mid d$ mais $p^2\nmid d$. On note $\widetilde{\mathfrak{M}}^{\lisse}_{2d}$ le compl\'ementaire de ces points.

\subsection{L'application des p\'eriodes}\label{parperiodes}

Notons $L_{d}$ le r\'eseau primitif\footnote{C'est le r\'eseau explicite $E_8(-1)^{\oplus 2}\oplus U^{\oplus 2}\oplus \langle -2d\rangle$, o\`u $U$ est le r\'eseau hyperbolique de rang\ $2$.} des surfaces K3 complexes munies d'une polarisation primitive de degr\'e $2d$, $G_d:=\SO(L_d)$, $G'_d:=\GSpin(L_d)$, et $\Gamma_d\subset G_d(\mathbb{Z})$ et $K_d\subset G_d(\widehat{\mathbb{Z}})$ les sous-groupes agissant trivialement sur le discriminant
$L_d^\vee/L_d$. Soit $$\Omega:=\{\omega\in \mathbb{P}(L_{d,\mathbb{C}})| \langle \omega,\omega\rangle=0,\langle\omega,\bar{\omega}\rangle>0 \}\simeq\SO(2,19)/\SO(2)\times \SO(19)$$
l'espace sym\'etrique\footnote{Au d\'etail pr\`es que $\Omega$ a deux composantes connexes.} associ\'e au groupe de Lie $\SO(2,19)$. 
L'\textbf{application des p\'eriodes}\footnote{Nous utilisons un rev\^etement double de l'application classique \cite{BeauvilleK3}.} pour les surfaces K3 est une immersion ouverte $P_{\mathbb{C}}:\widetilde{\mathfrak{M}}_{2d,\mathbb{C}}\to \Omega/\Gamma_d$ induite par $X\mapsto H^{2,0}(X)$.
On peut interpr\'eter le domaine de p\'eriodes $\Omega/\Gamma_d$ comme le complexifi\'e de la vari\'et\'e de Shimura
$\Sh_{d}:=\Sh_{K_d}(G_{d,\mathbb{Q}}, \Omega)$.

La donn\'ee de Shimura $\Omega$ se rel\`eve au groupe $G'_{d,\mathbb{Q}}$ \cite[4.2]{WeilK3}. On choisit comme niveau $K'_d\subset G'_d(\widehat{\mathbb{Z}})$ le sous-groupe de l'image r\'eciproque de $K_d$ dont la composante en $2$ appartient \`a un sous-groupe compact ouvert assez petit de $G'_d(\mathbb{Q}_2)$. Un tel choix permet d'assurer que $\Sh'_{d}:=\Sh_{K'_d}(G'_{d,\mathbb{Q}}, \Omega)$ est une vari\'et\'e quasi-projective.
On obtient un morphisme fini \'etale de vari\'et\'es de Shimura $\pi_{\mathbb{C}}:\Sh'_{d,\mathbb{C}}\to\Sh_{d,\mathbb{C}}$. 

Enfin, faisons agir $G'_{d,\mathbb{Q}}$ sur l'alg\`ebre de Clifford $H:=\Cl(L_{d,\mathbb{Q}})$ par multiplication \`a gauche. Il est possible de munir (non canoniquement) $H$ d'une forme symplectique $\psi$ $\Spin(L_{d,\mathbb{Q}})$-invariante, de sorte que l'immersion $G'_{\mathbb{Q}}\hookrightarrow\GSp(H,\psi)$ induise un morphisme fini et non ramifi\'e de vari\'et\'es de Shimura
$\iota_{\mathbb{C}} : \Sh'_{d,\mathbb{C}}\to\Ab_{\mathbb{C}}$, o\`u $\Ab$ est une vari\'et\'e de Shimura de type Siegel : un champ de modules de vari\'et\'es ab\'eliennes polaris\'ees. L'\textbf{application de Kuga-Satake} $\iota_{\mathbb{C}}$ induit un sch\'ema ab\'elien sur $\Sh'_{d,\mathbb{C}}$.

Notant $\widetilde{\mathfrak{M}}'_{2d,\mathbb{C}}:=\widetilde{\mathfrak{M}}_{2d,\mathbb{C}}\times_{\Sh_{d,\mathbb{C}}}\Sh'_{d,\mathbb{C}}$, on obtient le diagramme suivant :
\begin{equation*}
\xymatrix{
\widetilde{\mathfrak{M}}'_{2d,\mathbb{C}}\ar[d]\ar[r]^{P_{\mathbb{C}}'}
&  \Sh'_{d,\mathbb{C}}\ar[d]^{\pi_{\mathbb{C}}}\ar[r]^{\iota_{\mathbb{C}}}  &\Ab_{\mathbb{C}}  \\
\widetilde{\mathfrak{M}}_{2d,\mathbb{C}}\ar[r]^{P_{\mathbb{C}}}
& \Sh_{d,\mathbb{C}} .
}
\end{equation*}

Si $x\in \widetilde{\mathfrak{M}}'_{2d}(\mathbb{C})$ rel\`eve une surface K3 complexe polaris\'ee $[X]\in \mathfrak{M}_{2d}(\mathbb{C})$, la vari\'et\'e de Kuga-Satake de $X$ construite au paragraphe \ref{KSnaive} est la vari\'et\'e ab\'elienne $A$ correspondant \`a $\iota_{\mathbb{C}}(P'_{\mathbb{C}}(x))$.
De plus, le lien cohomologique (\ref{KSrel}) entre $X$ et $A$ se met en famille comme suit\footnote{Nous noterons toujours $\mathbf{H}$ le $H^1$ en famille de vari\'et\'es de Kuga-Satake, et $\mathbf{P}^2$ la cohomologie primitive en famille de surfaces K3 polaris\'ees. Nous noterons \'egalement $\mathbf{H}^{\otimes (1,1)}:=\mathbf{H}\otimes(\mathbf{H})^{\vee}=\End(\mathbf{H}).$}. Il existe une sous-variation de structures de Hodge naturelle $\mathbf{L}_{\B}\hookrightarrow\mathbf{H}^{\otimes(1,1)}_{\B}$ sur $\Sh'_{d,\mathbb{C}}$, et un isomorphisme $\mathbf{P}^2_{\B}(1)\xrightarrow{\sim}P_{\mathbb{C}}'^{*}\mathbf{L}_{\B}$ de variations de structures de Hodge sur $\widetilde{\mathfrak{M}}'_{2d,\mathbb{C}}$.

\`A l'aide des th\'eor\`emes de comparaison avec la cohomologie de de Rham alg\'ebrique (resp. la cohomologie $\ell$-adique), on obtient des isomorphismes compatibles
de fibr\'es vectoriels filtr\'es \`a connexion int\'egrable $\mathbf{P}^2_{\dR,\mathbb{C}}(1)\xrightarrow{\sim} P_{\mathbb{C}}'^*\mathbf{L}_{\dR,\mathbb{C}}\hookrightarrow P_{\mathbb{C}}'^*\mathbf{H}^{\otimes(1,1)}_{\dR,\mathbb{C}}$ (resp. de syst\`emes locaux $\ell$-adiques $\mathbf{P}^2_{\ell,\mathbb{C}}(1)\xrightarrow{\sim} P_{\mathbb{C}}'^*\mathbf{L}_{\ell,\mathbb{C}}\hookrightarrow P_{\mathbb{C}}'^*\mathbf{H}^{\otimes(1,1)}_{\ell,\mathbb{C}}$) sur $\widetilde{\mathfrak{M}}'_{2d,\mathbb{C}}$.

\subsection{Classes de Hodge absolues et rationalit\'e de la construction}\label{parQ}

Les corps reflex des trois vari\'et\'es de Shimura qui interviennent sont tous \'egaux \`a $\mathbb{Q}$ de sorte que celles-ci ont des mod\`eles canoniques $\Sh_d$, $\Sh'_d$ et $\Ab$ d\'efinis sur $\mathbb{Q}$, qui sont reli\'es par des morphismes de vari\'et\'es de Shimura $\pi_\mathbb{Q}$ et $\iota_{\mathbb{Q}}$.

L'outil essentiel utilis\'e pour descendre sur $\mathbb{Q}$ les constructions ci-dessus est la th\'eorie des classes de Hodge absolues, introduites implicitement \`a cet effet par Deligne \cite{WeilK3}. En effet, le morphisme (\ref{KSrel})
est induit par une classe de Hodge absolue.
Cela r\'esulte du th\'eor\`eme de Deligne que toute classe de Hodge sur une vari\'et\'e ab\'elienne est absolue \cite[I Theorem 2.11]{DMOS} ainsi que d'une r\'eduction au cas o\`u $X$ est une surface de Kummer, par le principe B de Deligne \cite[I Theorem 2.12]{DMOS}.

Ce fait permet d'abord de montrer \cite[2.2.1, 2.2.2]{Kisin} que les sous-fibr\'es vectoriels filtr\'es \`a connexion int\'egrable (resp. sous-syst\`emes locaux $\ell$-adiques) ci-dessus sont naturellement d\'efinis sur $\Sh'_{d}$ : on a $\mathbf{L}_{\dR,\mathbb{Q}}\hookrightarrow \mathbf{H}^{\otimes(1,1)}_{\dR,\mathbb{Q}}$ (resp. $\mathbf{L}_{\ell,\mathbb{Q}}\hookrightarrow \mathbf{H}^{\otimes(1,1)}_{\ell,\mathbb{Q}}$).

Il permet ensuite de v\'erifier \cite[Corollary 4.4]{MP} que l'application
des p\'eriodes $P_{\mathbb{C}}$ est d\'efinie sur $\mathbb{Q}$. Une preuve diff\'erente avait \'et\'e propos\'ee par Rizov \cite[Theorem 3.9.1]{RizovKS}, comme cons\'equence d'un th\'eor\`eme de la multiplication complexe pour les surfaces K3 de nombre de Picard $20$. Notant encore $\widetilde{\mathfrak{M}}'_{2d,\mathbb{Q}}:=\widetilde{\mathfrak{M}}_{2d,\mathbb{Q}}\times_{\Sh_{d}}\Sh'_{d}$, 
on obtient le diagramme :
$$\widetilde{\mathfrak{M}}'_{2d,\mathbb{Q}}\stackrel{P'_{\mathbb{Q}}}{\longrightarrow} \Sh'_{d}\stackrel{\iota_{\mathbb{Q}}}{\longrightarrow} \Ab.$$

Associ\'e \`a un argument de monodromie, il permet finalement \cite[Proposition 4.6]{MP} d'obtenir des isomorphismes
$\mathbf{P}^2_{\dR,\mathbb{Q}}(1)\xrightarrow{\sim} P_{\mathbb{Q}}'^*\mathbf{L}_{\dR,\mathbb{Q}}$ (resp. $\mathbf{P}^2_{\ell,\mathbb{Q}}(1)\xrightarrow{\sim} P_{\mathbb{Q}}'^*\mathbf{L}_{\ell,\mathbb{Q}}$) sur $\widetilde{\mathfrak{M}}'_{2d,\mathbb{Q}}$, qui sont les analogues de (\ref{KSrel}) en cohomologie de de Rham (resp. $\ell$-adique).

\subsection{Mod\`eles entiers de vari\'et\'es de Shimura}\label{parentier}

On souhaite \'etendre les constructions ci-dessus sur $\mathbb{Z}[\frac{1}{2}]$. Pour cela, il faut choisir des mod\`eles de nos vari\'et\'es sur $\mathbb{Z}[\frac{1}{2}]$.
Nous avons d\'ej\`a indiqu\'e que $\widetilde{\mathfrak{M}}_{2d,\mathbb{Q}}$ a un mod\`ele $\widetilde{\mathfrak{M}}_{2d}$ fourni par son interpr\'etation modulaire.
Il reste \`a construire des mod\`eles naturels pour les vari\'et\'es de Shimura qui interviennent. 

Une bonne notion, inspir\'ee par la d\'efinition des mod\`eles de N\'eron, a \'et\'e d\'egag\'ee par Milne \cite{Milne} (voir aussi \cite{Moonen}). La formuler n\'ecessite de ne pas travailler avec une structure de niveau fix\'ee. Pour cela, soient $(G,X)$ une donn\'ee de Shimura de corps reflex $E$, $\mathfrak{O}$ l'anneau des entiers de $E$, $p$ un nombre premier, $\mathfrak{p}$ une place de $E$ au-dessus de $p$, $\mathbb{A}_f^p$ les ad\`eles finis de $\mathbb{Q}$ dont la composante en $p$ est triviale, et $K_p\subset G(\mathbb{Q}_p)$ et $K^p\subset G(\mathbb{A}_f^p)$ des sous-groupes compacts ouverts. On consid\`ere la limite projective $\Sh_{K_p}(G,X):=\varprojlim_{K^p}\Sh_{K_pK^p}(G,X)$ sur les structures de niveau premi\`eres \`a $p$ : c'est un $E$-sch\'ema muni d'une action continue de $G(\mathbb{A}_f^p)$. C'est pour ce sch\'ema qu'on peut esp\'erer l'existence d'un mod\`ele v\'erifiant une propri\'et\'e universelle de type N\'eron.

\begin{defi}\label{defint}
Un \textbf{mod\`ele entier canonique lisse} en $\mathfrak{p}$ de $\Sh_{K_p}(G,X)$ est un $\mathfrak{O}_{(\mathfrak{p})}$-sch\'ema s\'epar\'e, r\'egulier et formellement lisse $\mathpzc{Sh}_{K_p}(G,X)$ muni d'une action continue de $G(\mathbb{A}_f^p)$ et d'une identification $G(\mathbb{A}_f^p)$-\'equivariante de sa fibre g\'en\'erique avec $\Sh_{K_p}(G,X)$ tel que pour tout 
$\mathfrak{O}_{(\mathfrak{p})}$-sch\'ema r\'egulier et formellement lisse $\mathcal{X}$, tout morphisme $\mathcal{X}_E\to\Sh_{K_p}(G,X)$ se prolonge en un morphisme $\mathcal{X}\to \mathpzc{Sh}_{K_p}(G,X)$.

Si un tel mod\`ele existe, il est unique. On d\'efinit alors le mod\`ele entier canonique lisse de $\Sh_{K_pK^p}(G,X)$ : c'est $\mathpzc{Sh}_{K_pK^p}(G,X):=\mathpzc{Sh}_{K_p}(G,X)/K^p$.
\end{defi}

Dans les cas qui nous int\'eressent, de tels mod\`eles ont \'et\'e construits par Kisin \cite{Kisin} quand $p\nmid 2d$, et par Madapusi Pera \cite{MP2} quand $p\neq 2$ :

\begin{theo}
Les vari\'et\'es de Shimura $\Sh_d$ et $\Sh'_d$ admettent des mod\`eles $\mathpzc{Sh}_d$ et $\mathpzc{Sh}'_d$ sur $\mathbb{Z}[\frac{1}{2}]$ qui sont canoniques lisses sur $\mathbb{Z}_{(p)}$ pour tout $p\neq 2$.
\end{theo}

Si $p\nmid d$, le principe de la preuve de ce th\'eor\`eme remonte \`a Milne \cite{Milne}, et une telle preuve a \'et\'e annonc\'ee par Vasiu \cite{Vasiu}. La vari\'et\'e de Shimura $\Ab$ a une interpr\'etation modulaire qui permet de construire un mod\`ele entier
naturel $\mathpzc{Ab}$ de $\Ab$ sur $\mathbb{Z}[\frac{1}{2}]$ : $\mathpzc{Ab}$ est encore un champ de modules de vari\'et\'es ab\'eliennes polaris\'ees.
On d\'efinit alors $\mathpzc{Sh}'_d$ comme la normalisation de l'adh\'erence de $\Sh'_d$ dans $\mathpzc{Ab}$\footnote{Nous sommes ici impr\'ecis sur le choix des structures de niveau.\label{fn:repeat}}. La difficult\'e, surmont\'ee par Kisin, est de v\'erifier que $\mathpzc{Sh}'_d$ est bien lisse. La propri\'et\'e universelle des mod\`eles canoniques lisses r\'esulte alors de th\'eor\`emes d'extension de sch\'emas ab\'eliens. Kisin construit enfin $\mathpzc{Sh}_d$ \`a partir de $\mathpzc{Sh}'_d$ de sorte que $\pi:\mathpzc{Sh}'_d\to \mathpzc{Sh}_d$ soit toujours fini \'etale, en montrant que les difficult\'es identifi\'ees par Moonen \cite[3.21]{Moonen} n'apparaissent pas.

Pour r\'eduire le cas g\'en\'eral \`a cette situation, Madapusi Pera plonge primitivement le r\'eseau $L_d\otimes \mathbb{Z}_{(p)}$ dans un $\mathbb{Z}_{(p)}$-r\'eseau $M$ de discriminant premier \`a $p$ \cite[Lemma 6.8]{MP2}.
Les vari\'et\'es de Shimura associ\'ees \`a $\SO(M_{\mathbb{Q}})$ et $\GSpin(M_{\mathbb{Q}})$ sont justiciables des r\'esultats de Kisin et poss\`edent donc des mod\`eles entiers canoniques lisses en $p$. Madapusi Pera est capable de d\'ecrire avec assez de pr\'ecision l'adh\'erence de $\Sh_d$ et $\Sh'_d$ dans ces mod\`eles pour en d\'eduire l'existence des mod\`eles entiers canoniques lisses $\mathpzc{Sh}_d$ et $\mathpzc{Sh}'_d$.
A~posteriori, la construction montre qu'on aurait pu d\'efinir $\mathpzc{Sh}'_d$ comme le lieu lisse de la normalisation de l'adh\'erence $\Sh'_d$ dans $\mathpzc{Ab}$
\footref{fn:repeat}.

\medskip

La propri\'et\'e universelle des mod\`eles canoniques lisses, appliqu\'ee \`a une limite projective $\mathcal{X}$ d'espaces de modules de surfaces K3 avec structures de niveau, permet d'\'etendre l'application des p\'eriodes $P_{\mathbb{Q}}$ en $P: \widetilde{\mathfrak{M}}_{2d}^{\lisse}\to \mathpzc{Sh}_{d}$ \cite[Proposition 4.7]{MP}. Notant toujours $\widetilde{\mathfrak{M}}'^{\lisse}_{2d}:=\widetilde{\mathfrak{M}}^{\lisse}_{2d}\times_{\mathpzc{Sh}_{d}}\mathpzc{Sh}'_{d}$, on obtient le diagramme :
$$ \widetilde{\mathfrak{M}}_{2d}'^{\lisse}\stackrel{P'}{\longrightarrow}\mathpzc{Sh}'_{d}\stackrel{\iota}{\longrightarrow} \mathpzc{Ab}.$$

  Il est facile d'\'etendre $\mathbf{L}_{\ell,\mathbb{Q}}\hookrightarrow \mathbf{H}^{\otimes(1,1)}_{\ell,\mathbb{Q}}$ en un sous-syst\`eme local $\ell$-adique $\mathbf{L}_{\ell}\hookrightarrow\mathbf{H}^{\otimes(1,1)}_{\ell}$ sur $\mathpzc{Sh}'_{d,\mathbb{Z}[\frac{1}{2\ell}]}$ satisfaisant $\mathbf{P}^2_{\ell}(1)\xrightarrow{\sim} P'^*\mathbf{L}_{\ell}$. La construction de Kisin et Madapusi Pera montre qu'il est \'egalement possible d'\'etendre naturellement $\mathbf{L}_{\dR,\mathbb{Q}}\hookrightarrow \mathbf{H}^{\otimes(1,1)}_{\dR,\mathbb{Q}}$ en un sous-fibr\'e vectoriel filtr\'e \`a connexion int\'egrable $\mathbf{L}_{\dR}\hookrightarrow\mathbf{H}^{\otimes(1,1)}_{\dR}$ sur $ \mathpzc{Sh}'_{d}$. La compatibilit\'e $\mathbf{P}^2_{\dR}(1)\xrightarrow{\sim} P'^*\mathbf{L}_{\dR}$ sera v\'erifi\'ee \`a la proposition \ref{cristaKS}.

\begin{rema}\label{remspec}
Quand $p\mid d$ mais $p^2\nmid d$, Madapusi Pera \cite{MP2} consid\`ere une variante des d\'efinitions et constructions ci-dessus et obtient des mod\`eles entiers canoniques en $p$ de $\Sh_d$ et $\Sh'_d$ non n\'ecessairement lisses,
qui permettent d'\'etendre l'application des p\'eriodes \`a $\widetilde{\mathfrak{M}}_{2d}$ tout entier \cite[Proposition 4.7]{MP}.
\end{rema}

\subsection{Le th\'eor\`eme de Torelli}

Venons-en \`a pr\'esent au th\'eor\`eme principal de cette partie, qui g\'en\'eralise le th\'eor\`eme de Torelli en caract\'eristique mixte :

\begin{theo}
\label{Torelli}
L'application des p\'eriodes $P$ est une immersion ouverte.
\end{theo}

Cet \'enonc\'e avait d\'ej\`a \'et\'e obtenu par Maulik \cite[Proposition 5.10]{Maulik} quand $p\nmid d$ et $p\geq 5$. Le cas g\'en\'eral, d\^u \`a Madapusi Pera \cite[Theorem 4.8]{MP}, suit une strat\'egie proche de celle de Maulik. Le choix d'une construction de Kuga-Satake utilisant toute l'alg\`ebre de Clifford est une simplification technique importante introduite dans \cite[Proposition 13]{Charles}.
Signalons \'egalement un r\'esultat ant\'erieur, moins pr\'ecis, de Rizov \cite[Corollary 7.2.3]{RizovKS}.

La preuve suit les arguments du th\'eor\`eme de Torelli infinit\'esimal classique. Les d\'eformations infinit\'esimales d'une vari\'et\'e lisse $X$ sont contr\^ol\'ees par le groupe $H^1(X,T_X)$, \`a l'aide de l'application de Kodaira-Spencer. Si $X$ est une surface K3, $H^1(X, T_X)\simeq H^1(X,\Omega^1_X)$ est un gradu\'e pour la filtration de Hodge de la cohomologie de de Rham alg\'ebrique de $X$. Ceci explique que la preuve du th\'eor\`eme \ref{Torelli} soit une cons\'equence d'un r\'esultat de compatibilit\'e en cohomologie de Rham de la construction de Kuga-Satake :

\begin{prop}\label{cristaKS}
L'isomorphisme $\mathbf{P}^2_{\dR,\mathbb{Q}}(1)\xrightarrow{\sim} P_{\mathbb{Q}}'^*\mathbf{L}_{\dR,\mathbb{Q}}$ induit un isomorphisme 
$\mathbf{P}^2_{\dR}(1)\xrightarrow{\sim} P'^*\mathbf{L}_{\dR}$ de fibr\'es vectoriels filtr\'es \`a connexion int\'egrable.
\end{prop}

\noindent{\sc Preuve du th\'eor\`eme \ref{Torelli}} ---
Soient $[X]\in \widetilde{\mathfrak{M}}_{2d}'^{\lisse}(K)$ un point correspondant \`a une surface K3 polaris\'ee $(X,\xi)$ sur un corps alg\'ebriquement clos $K$ de caract\'eristique $\geq 3$, et $v\in \widetilde{\mathfrak{M}}_{2d}'^{\lisse}(K[\varepsilon]/\varepsilon^2)$ un vecteur tangent \`a $\widetilde{\mathfrak{M}}_{2d}'^{\lisse}$ en $[X]$ contract\'e par l'application des p\'eriodes $P'$. Comme $\mathbf{P}^2_{\dR}(1)\simeq P'^*\mathbf{L}_{\dR}$ par la proposition \ref{cristaKS}, l'application $\nabla_v:\mathbf{P}^2_{\dR}|_{[X]}\to \mathbf{P}^2_{\dR}|_{[X]}$ est nulle. L'application induite par transversalit\'e de Griffiths entre les gradu\'es pour la filtration de Hodge $\nabla_v:H^0(X,\Omega_X^2)\to H^1(X,\Omega^1_X)^{\perp\xi}$ l'est donc \'egalement. Comme celle-ci est donn\'ee par le cup-produit avec la classe de Kodaira-Spencer $\textsc{ks}(v)\in H^1(X,T_X)$, et que $\Omega_X^2\simeq \mathcal{O}_X$, il s'ensuit que cette classe de Kodaira-Spencer est nulle. Par propri\'et\'e de l'espace de modules, $v$ est donc nul.

Comme $\widetilde{\mathfrak{M}}_{2d}'^{\lisse}$ et $\mathpzc{Sh}'_{d}$ sont lisses de la m\^eme dimension relative $19$ sur $\Spec(\mathbb{Z}[\frac{1}{2}])$, nous avons montr\'e que $P'$ est \'etale, donc que $P$ est \'etale. Comme $\widetilde{\mathfrak{M}}_{2d}^{\lisse}$ est s\'epar\'e et que $P|_{\widetilde{\mathfrak{M}}_{2d,\mathbb{Q}}^{\lisse}}$ est une immersion ouverte par le th\'eor\`eme de Torelli en caract\'eristique nulle, le Main Theorem de Zariski montre que $P$ est une immersion ouverte.
 \qed

\medskip

La preuve de la proposition \ref{cristaKS} exploite l'interpr\'etation cristalline de la cohomologie de de Rham, en se reposant sur des versions enti\`eres du th\'eor\`eme de comparaison entre cohomologie \'etale $p$-adique et cohomologie cristalline. \`A cet effet, Maulik et Charles s'appuyaient sur la th\'eorie de Fontaine-Laffaille et Fontaine-Messing \cite{FontaineMessing}, et sur des am\'eliorations de celle-ci dues \`a Kisin \cite{Kisinbetter}. Cela ne leur permettait pas de traiter les caract\'eristiques $2$ et $3$. Madapusi Pera \'evite \`a cette \'etape toute hypoth\`ese sur la caract\'eristique, en remarquant qu'il n'a besoin de consid\'erer que des vari\'et\'es ordinaires, et en exploitant un r\'esultat ant\'erieur plus faible de Bloch et Kato \cite{BlochKato}, au sujet duquel on pourra aussi consulter l'appendice de \cite{PR}.

\smallskip

\noindent{\sc Preuve de la proposition \ref{cristaKS}} ---  
Soient $K$ un corps alg\'ebriquement clos de carac\-t\'eristique $p\geq 3$ et $s\in\widetilde{\mathfrak{M}}_{2d}'^{\lisse}(W(K))$ dont la fibre sp\'eciale correspond \`a une surface~K3~$X$ sur $K$. Les groupes de cohomologie de de Rham $\mathbf{P}^2_{\dR}|_{s}$ (resp. $P'^*\mathbf{L}_{\dR}(-1)|_{s}$) sont des $W(K)$-modules qui ont une structure de $F$-cristaux\footnote{Nous avons modifi\'e les twists \`a la Tate car $\mathbf{P}^2_{\dR}(1)|_s$, qui peut avoir des pentes n\'egatives, n'est pas vraiment un $F$-cristal.} (i.e. qui portent une action semi-lin\'eaire du Frobenius) induite par leur interpr\'etation cristalline. Ogus \cite[\S 7]{OgusII} a montr\'e que l'isomorphisme $\mathbf{P}^2_{\dR}|_{s_\mathbb{Q}}\xrightarrow{\sim} P_{\mathbb{Q}}'^*\mathbf{L}_{\dR}(-1)|_{s_\mathbb{Q}}$ commute \`a cette action du Frobenius. Comme expliqu\'e par Madapusi Pera \cite[Lemma 4.9]{MP}, on peut aujourd'hui voir cet \'enonc\'e comme une cons\'equence d'un th\'eor\`eme de Blasius et Wintenberger (\cite{Blasius}, \cite[Theorem 5.6.3]{Moonen}) selon lequel le th\'eor\`eme de comparaison entre cohomologie cristalline et cohomologie \'etale $p$-adique est compatible aux r\'ealisations des classes de Hodge absolues sur les vari\'et\'es ab\'eliennes.



Supposons $X$ ordinaire. Les polygones de Hodge et de Newton du $F$-cristal $\mathbf{P}^2_{\dR}|_{s}$ co\"incident alors. Comme la vari\'et\'e de Kuga-Satake $A$ associ\'ee \`a $X$ est \'egalement ordinaire (\cite[Proposition 2.5]{Nygaard}, \cite[Theorem 7.8]{OgusII}) et que  $P'^*\mathbf{L}_{\dR}(-1)|_s$ est un sous-$F$-cristal primitif de $H^{\otimes (1,1)}_{\cris}(A)(-1)$, il en va de m\^eme pour le $F$-cristal $P'^*\mathbf{L}_{\dR}(-1)|_{s}$.
La d\'ecomposition de Newton-Hodge \cite[Theorem 1.6.1]{Katz} montre alors que le $F$-cristal $\mathbf{P}^2_{\dR}|_{s}$ (resp. $P'^*\mathbf{L}_{\dR}(-1)|_s$) est canoniquement somme directe de trois sous-$F$-cristaux de pentes $0$, $1$ et $2$. Un th\'eor\`eme de Bloch et Kato \cite[Theorem 9.6.2]{BlochKato} montre que 
ces sous-$F$-cristaux s'identifient aux gradu\'es tensoris\'es par $W(K)$ d'une filtration naturelle sur $\mathbf{P}^2_{p}|_{s_\mathbb{Q}}$ (resp. $P_{\mathbb{Q}}'^*\mathbf{L}_{p}(-1)|_{s_\mathbb{Q}}$). Comme $\mathbf{P}^2_{p}|_{s_\mathbb{Q}}\xrightarrow{\sim} P_{\mathbb{Q}}'^*\mathbf{L}_{p}(-1)|_{s_\mathbb{Q}}$, on conclut que $\mathbf{P}^2_{\dR}\to P'^*\mathbf{L}_{\dR}(-1)$ est un isomorphisme. 

Appliquant ce r\'esultat \`a des rel\`evements des points g\'en\'eriques g\'eom\'etriques de la fibre sp\'eciale de $\widetilde{\mathfrak{M}}_{2d}'^{\lisse}$, on voit que l'isomorphisme
$\mathbf{P}^2_{\dR,\mathbb{Q}}(1)\xrightarrow{\sim} P_{\mathbb{Q}}'^*\mathbf{L}_{\dR,\mathbb{Q}}$ se prolonge en un isomorphisme de fibr\'es vectoriels sur un ouvert contenant les points g\'en\'eriques de la fibre sp\'eciale de $\widetilde{\mathfrak{M}}_{2d}'^{\lisse}$. Comme $\widetilde{\mathfrak{M}}_{2d}'^{\lisse}$ est normal, il s'\'etend en un isomorphisme sur $\widetilde{\mathfrak{M}}_{2d}'^{\lisse}$ tout entier. Cet isomorphisme pr\'eserve la connexion et la filtration, car c'est le cas sur la fibre g\'en\'erique.
\qed
\begin{rema}
La preuve pr\'ec\'edente montre en fait \cite[Corollary 4.14]{MP} que l'on a un isomorphisme $\mathbf{P}^2_{\cris}\xrightarrow{\sim} P'^*\mathbf{L}_{\cris}(-1)\hookrightarrow \mathbf{H}^{\otimes(1,1)}_{\cris}(-1)$ de $F$-cristaux sur $\widetilde{\mathfrak{M}}_{2d,\mathbb{F}_p}'^{\lisse}$, qui met en famille l'analogue cristallin de la relation (\ref{KSrel}).
\end{rema}

\begin{rema}\label{global}
Le th\'eor\`eme \ref{Torelli} implique que $\widetilde{\mathfrak{M}}_{2d,\mathbb{F}_p}^{\lisse}$ a pour espace de modules grossier un sch\'ema quasi-projectif. Il permet \'egalement de montrer que celui-ci est g\'eom\'etriquement int\`egre si $p^2\nmid d$ \cite[Corollary 4.16]{MP}.
\end{rema}

\begin{rema}
On ne peut esp\'erer que l'application des p\'eriodes $P$ soit surjective : il faudrait au moins prendre en compte les surfaces K3 quasi-polaris\'ees, c'est-\`a-dire munies d'un fibr\'e en droites gros et nef. Les articles \cite{Maulik} et \cite{MP} se placent dans ce cadre l\'eg\`erement plus g\'en\'eral. 
La question de la surjectivit\'e de l'application des p\'eriodes est alors intimement li\'ee \`a l'existence d'un crit\`ere de bonne r\'eduction potentielle du type N\'eron-Ogg-Shafarevich pour les surfaces K3. 
Matsumoto \cite{Matsumoto} a obtenu un tel r\'esultat en caract\'eristique assez grande, en s'inspirant de techniques de Maulik \cite[\S 4]{Maulik}.
\end{rema}

\begin{rema}
Un th\'eor\`eme de Torelli \`a l'\'enonc\'e plus proche des formulations classiques (la surface K3 est d\'etermin\'ee par des donn\'ees cohomologiques) avait \'et\'e obtenu par Ogus \cite{O2} pour les surfaces K3 Shioda-supersinguli\`eres en caract\'eristique $\geq 5$.
\end{rema}

\section{Construction de fibr\'es en droites}

Dans cette partie, nous expliquons les id\'ees de la preuve de la conjecture de Tate pour les surfaces K3 en caract\'eristique diff\'erente de $2$ :

\begin{theo}\label{thTate2}
Soit $X$ une surface K3 sur un corps de type fini $k$ de carac\-t\'eristique~$\neq 2$. Soient $\ell$ un nombre premier inversible dans $k$, $\bar{k}$ une cl\^oture s\'eparable de $k$ et \mbox{$\Gamma_k:=\Gal(\bar{k}/k)$.} Alors l'application classe de cycle induit un isomorphisme :
$$\Pic(X)\otimes \mathbb{Q}_\ell\to H^2_{\et}(X_{\bar{k}},\mathbb{Q}_\ell(1))^{\Gamma_k}.$$
\end{theo}

Les premiers r\'esultats positifs ont \'et\'e obtenus pour des surfaces K3 elliptiques (\cite[Theorem 5.2]{ASD}, \cite[Theorem 1.7]{Artin}) ou de degr\'e $2$ \cite[Theorem 4]{RSZ}. 

Toutes les avanc\'ees ult\'erieures ont repos\'e sur la construction de Kuga-Satake. En caract\'eristique nulle, celle-ci permet facilement de se ramener \`a la conjecture de Tate pour les vari\'et\'es ab\'eliennes prouv\'ee par Faltings \cite{Faltings,Faltingscomplements}. Sur les corps finis, Nygaard \cite{Nygaard} (resp. Nygaard et Ogus \cite{NO}) ont \'etendu ces arguments, et ramen\'e la conjecture de Tate pour les surfaces K3 ordinaires (resp. de hauteur finie en carac\-t\'eristique $\geq 5$)  \`a la conjecture de Tate pour les vari\'et\'es ab\'eliennes sur les corps finis prouv\'ee par Tate \cite{Tateendo}. Apr\`es un premier paragraphe de g\'en\'eralit\'es, nous rappelons bri\`evement ces travaux aux paragraphes \ref{Tnulle} et \ref{Tcan}.

\vspace{1em}

C'est Maulik \cite{Maulik} qui a le premier obtenu des r\'esultats dans le cas supersingulier, plus difficile, et laiss\'e ouvert par Nygaard et Ogus. Plus pr\'ecis\'ement, Maulik montre le th\'eor\`eme \ref{thTate2} quand $k$ est un corps fini de caract\'eristique $p\geq 5$, $X$ est supersinguli\`ere et $X$ poss\`ede une polarisation de degr\'e $2d$ avec $p>2d+4$. Sa strat\'egie est sp\'ecifique aux surfaces K3 supersinguli\`eres, et proc\`ede par r\'eduction au cas elliptique. L'originalit\'e de sa m\'ethode r\'eside dans l'utilisation formes automorphes construites \`a l'aide de travaux de Borcherds \cite{Borcherds}.

Dans \cite{Charles}, Charles a poursuivi et am\'elior\'e cette strat\'egie, en exploitant davantage la construction de Kuga-Satake. Cela lui permet de contourner les difficult\'es li\'ees \`a la mauvaise r\'eduction de surfaces K3 rencontr\'ees par Maulik. Nous expliquons comment les arguments qu'il d\'eveloppe
  permettent de montrer le th\'eor\`eme \ref{thTate2} quand $k$ est un corps fini de caract\'eristique $ p\geq 5$, $X$ est supersinguli\`ere et $X$ poss\`ede une polarisation de degr\'e premier \`a $p$
\footnote{L'article \cite{Charles} a pour th\'eor\`eme principal la conjecture de Tate pour les surfaces K3 supersinguli\`eres sur un corps fini de caract\'eristique $\geq 5$. Il contient malheureusement une erreur. Dans la preuve de \cite[Proposition 25]{Charles}, il faut v\'erifier que $\overline{T}$ est lisse pour pouvoir appliquer \cite[Lemma 5.12]{Maulik}. Sous les hypoth\`eses de \cite[Proposition 25]{Charles}, c'est vrai, et d\^u \`a Kisin \cite{Kisin}. Cependant, dans la preuve de \cite[Corollary 29]{Charles}, la proposition \cite[Proposition 25]{Charles} est utilis\'ee dans une g\'en\'eralit\'e o\`u \cite{Kisin} ne s'applique pas.

  Un erratum \`a \cite{Charles} expliquera comment corriger ces arguments. Dans ce texte, nous nous contentons de traiter le cas o\`u il existe une polarisation de degr\'e premier \`a la caract\'eristique.}.

Le paragraphe \ref{MauCha}, dont le r\'esultat principal est le th\'eor\`eme \ref{thTate3}, est consacr\'e \`a cette approche.

\vspace{1em}

 Dans les paragraphes \ref{Tspe} et \ref{TMP}, nous discutons la preuve par Madapusi Pera \cite{MP} du th\'eor\`eme \ref{thTate2}. Sa strat\'egie est tout \`a fait dans le prolongement de la preuve de la conjecture de Tate en caract\'eristique nulle, ainsi que des travaux de Nygaard et Ogus. Elle a l'avantage de ne pas distinguer entre surfaces K3 supersinguli\`eres et surfaces K3 de hauteur finie, de ne pas n\'ecessiter une r\'eduction au cas elliptique, et de permettre de traiter le cas des corps de type fini.
Elle s'appuie de mani\`ere essentielle sur les travaux r\'ecents de Kisin \cite{Kisinnew} sur la r\'eduction modulo $p$ de vari\'et\'es de Shimura de type Hodge.

\subsection{Pr\'eliminaires}

 Pour prouver le th\'eor\`eme \ref{thTate2}, on peut remplacer $k$ par une extension finie $k'$. En effet, si la conjecture de Tate vaut pour $X_{k'}$ et si $\alpha\in H^2_{\et}(X_{\bar{k}},\mathbb{Q}_\ell(1))^{\Gamma_k}$, on peut \'ecrire  $\alpha$ comme $\mathbb{Q}_\ell$-combinaison lin\'eaire $\sum_i \lambda_i[\mathcal{L}_i]_{\ell}$ de classes de fibr\'e en droites sur $X_{k'}$. On a alors $[k':k]\alpha= \sum_i \lambda_i[N_{k'/k}(\mathcal{L}_i)]_{\ell}$.

Pour cette raison, on s'autorisera dans la suite \`a remplacer $k$ par une extension finie sans n\'ecessairement le mentionner explicitement. On choisit en particulier une telle extension de sorte que tous les fibr\'es en droites sur $X_{\bar{k}}$ soient d\'efinis sur $k$, et on munit $X$ d'une polarisation primitive $\xi$ de degr\'e $2d$.

\medskip

Nous avions exclu au paragraphe \ref{edm} un nombre fini de surfaces K3 supersp\'eciales. Ce n'est pas grave pour deux raisons. D'une part, ce n'\'etait pas vraiment n\'ecessaire (voir la remarque \ref{remspec}). D'autre part, celles-ci sont supersinguli\`eres et poss\`edent des d\'eformations supersinguli\`eres non supersp\'eciales \cite[Remark 2.7]{O1}. Par \cite[Theorem 1.1]{Artin}, il suffit de prouver la conjecture de Tate pour une telle d\'eformation. On peut donc supposer $[X]\in\mathfrak{M}_{2d}^{\lisse}(k)$.

\medskip

Rappelons le diagramme $\widetilde{\mathfrak{M}}_{2d}'^{\lisse}\stackrel{P'}{\longrightarrow}\mathpzc{Sh}'_{d}\stackrel{\iota}{\longrightarrow} \mathpzc{Ab}$ obtenu au paragraphe \ref{parentier}. L'application des p\'eriodes $P'$ est une immersion ouverte par le th\'eor\`eme de Torelli \ref{Torelli}, et, par construction, l'application de Kuga-Satake $\iota$ est quasi-finie, et finie si $p\nmid d$.

Quitte \`a remplacer $k$ par une extension finie, on peut supposer que $[X]\in\mathfrak{M}_{2d}^{\lisse}(k)$ se rel\`eve en un point encore not\'e $[X]\in\widetilde{\mathfrak{M}}_{2d}'^{\lisse}(k)$. On dira que la vari\'et\'e ab\'elienne $A$ correspondant \`a $\iota(P'([X]))$ est la vari\'et\'e de Kuga-Satake
associ\'ee \`a $X$.

Explicitons les informations cohomologiques obtenues dans la partie pr\'ec\'edente. Si $x\in \mathpzc{Sh}'_{d}(k)$, et si $A$ est la vari\'et\'e ab\'elienne associ\'ee \`a $\iota(x)$, on dispose, pour diverses th\'eories cohomologiques, d'un sous-espace canonique $L(A)\hookrightarrow \End(H^1(A))$. Si de plus $x=P'([X])$, on a construit, pour ces m\^emes th\'eories cohomologiques, un isomorphisme $P^2(X)(1)\xrightarrow{\sim}L(A)$.
En particulier, si $A$ est la vari\'et\'e de Kuga-Satake associ\'ee \`a $X$, on dispose d'une inclusion naturelle
$P^2(X)(1)\hookrightarrow \End(H^1(A))$.


\subsection{La conjecture de Tate en caract\'eristique nulle}\label{Tnulle}

Expliquons maintenant la preuve classique du th\'eor\`eme \ref{thTate2} en caract\'eristique nulle  (\cite[Theorem 5.6 (a)]{Tatemotives}, \cite[Theorem 1.6.1 (ii)]{Andre}). Fixant un plongement $\bar{k}\hookrightarrow \mathbb{C}$, on obtient
un diagramme commutatif :
\begin{equation*}
\xymatrix{
P^2_{\B}(X_{\mathbb{C}},\mathbb{Q}(1))\ar@{^{(}->}[r]^{j_{\B}\hspace{1.3em}}\ar[d]_{}
&  \End(H^1_{\B}(A_{\mathbb{C}},\mathbb{Q})) \ar[d]_{}   \\
P^2_{\et}(X_{\bar{k}},\mathbb{Q}_{\ell}(1))\ar@{^{(}->}[r]^{j_{\ell}\hspace{1.3em}}
& \End(H^1_{\et}(A_{\bar{k}},\mathbb{Q}_{\ell})).
}
\end{equation*}
La ligne sup\'erieure est un morphisme de structures de Hodge rationnelles, la ligne inf\'erieure est $\Gamma_k$-\'equivariante, et les fl\`eches verticales sont donn\'ees par le th\'eor\`eme de comparaison entre cohomologie de Betti et cohomologie $\ell$-adique. Comme la structure de Hodge $ \End(H^1_{\B}(A_{\mathbb{C}},\mathbb{Q}))$ est polarisable, il existe une r\'etraction $q_{\B}: \End(H^1_{\B}(A_{\mathbb{C}},\mathbb{Q}))\to P^2_{\B}(X_{\mathbb{C}},\mathbb{Q}(1))$ de $j_{\B}$ qui est un morphisme de structures de Hodge. On note $q_{\ell}$ son extension des scalaires \`a $\mathbb{Q}_{\ell}$.

Soit $\alpha\in H^2_{\et}(X_{\bar{k}},\mathbb{Q}_{\ell}(1))^{\Gamma_k}$. Quitte \`a la modifier par un multiple de la polarisation, on peut supposer $\alpha\in P^2_{\et}(X_{\bar{k}},\mathbb{Q}_{\ell}(1))^{\Gamma_k}$. Son image $j_\ell(\alpha)$ est encore $\Gamma_k$-invariante. Par la conjecture de Tate pour les vari\'et\'es ab\'eliennes prouv\'ee par Faltings \cite{Faltings,Faltingscomplements},
cette image est $\mathbb{Q}_{\ell}$-combinaison lin\'eaire de classes d'endomorphismes $f_i$ de $A$ : \'ecrivons $j_\ell(\alpha)=\sum_i \lambda_i[f_i]_{\ell}$. La classe $q_B([f_i]_{\B})$ est de Hodge, et par le th\'eor\`eme des classes $(1,1)$ de Lefschetz, c'est la classe d'un fibr\'e en droites $\mathcal{L}_i$ sur $X$. On calcule alors $\alpha=q_{\ell}(j_{\ell}(\alpha))=q_{\ell}(\sum_i \lambda_i[f_i]_{\ell})=\sum_i \lambda_iq_{\ell}([f_i]_{\ell})=\sum_i\lambda_i[\mathcal{L}_i]_{\ell}$.

\subsection{Rel\`evements canoniques et quasi-canoniques}\label{Tcan}

Les articles \cite{Nygaard} et \cite{NO} adaptent cette strat\'egie. Supposons $k$ fini et $X$ ordinaire, et esquissons les arguments de \cite{Nygaard}.
   Nygaard (voir aussi \cite{canonique}) montre l'existence d'un rel\`evement polaris\'e $\mathcal{X}$ de $X$ sur $W(k)$, dit canonique, qui v\'erifie la propri\'et\'e suivante. Soit $\mathcal{A}$ le sch\'ema ab\'elien de Kuga-Satake associ\'e \`a $\mathcal{X}$ : alors $A$ est ordinaire et $\mathcal{A}$ est son rel\`evement canonique\footnote{Nygaard \cite[Corollary 2.5, Proposition 2.8]{Nygaard} montre ce r\'esultat \`a isog\'enie pr\`es, mais voir \cite[Proposition 7.2.2]{RizovKS} et \cite[Proposition 4.22]{MP}.}. En particulier,  tout endomorphisme de $A$ se rel\`eve en un endomorphisme de $\mathcal{A}$.
   
   La preuve du paragraphe pr\'ec\'edent fonctionne alors : une classe dans $P^2_{\et}(X_{\bar{k}},\mathbb{Q}_{\ell}(1))^{\Gamma_k}$ induit une classe dans $\End(H^1_{\et}(A_{\bar{k}},\mathbb{Q}_{\ell}))^{\Gamma_k}$. Par la conjecture de Tate pour les vari\'et\'es ab\'eliennes prouv\'ee par Tate \cite{Tateendo}, celle-ci est $\mathbb{Q}_{\ell}$-combinaison lin\'eaire de classes d'endomorphismes de $A$. Ces endomorphismes se rel\`event \`a $\mathcal{A}$, o\`u leurs classes de Betti sont de Hodge, et induisent par projection des classes de Hodge dans $P^2_{\B}(\mathcal{X}_{\mathbb{C}},\mathbb{Q}(1))$. Par le th\'eor\`eme des classes $(1,1)$ de Lefschetz, celles-ci correspondent \`a des fibr\'es en droites qu'on peut sp\'ecialiser sur $X$ pour conclure.
   
La strat\'egie de \cite{NO} est similaire mais plus compliqu\'ee. Ils travaillent avec la cohomologie cristalline plut\^ot qu'avec la cohomologie $\ell$-adique, et d\'eveloppent, pour les surfaces K3 de hauteur finie, une notion plus faible de rel\`evement quasi-canonique. Comme les rel\`evements canoniques de \cite{Nygaard}, ces rel\`evements quasi-canoniques ont la propri\'et\'e remarquable que tout fibr\'e en droites sur $X$ s'y rel\`eve.

\subsection{Endomorphismes sp\'eciaux}\label{Tspe}

La difficult\'e pour appliquer ces id\'ees \`a une surface K3 supersinguli\`ere est qu'il n'est pas possible de trouver un rel\`evement de celle-ci en caract\'eristique nulle auquel tous les fibr\'es en droites se rel\`event.

On peut tout de m\^eme exploiter les arguments ci-dessus en consid\'erant plusieurs rel\`evements distincts de $X$. Cela permet  de relier les fibr\'es en droites sur $X$ \`a certains endomorphismes de $A$, dits sp\'eciaux. Nous suivons la pr\'esentation syst\'ematique de Madapusi Pera \cite{MP2,MP} ; ces arguments sont \'egalement utilis\'es dans la preuve de \cite[Proposition 22]{Charles}.

\begin{defi}
Soient $k$ un corps et $A$ la vari\'et\'e ab\'elienne associ\'ee \`a un $k$-point de $\mathpzc{Sh}'_d$. Un endomorphisme $f$ de $A$ est dit \textbf{sp\'ecial} si sa classe appartient \`a $L(A)\hookrightarrow\End(H^1(A))$ pour toutes les cohomologies $\ell$-adiques avec $\ell$ inversible dans $k$, et en cohomologie cristalline si $k$ est de caract\'eristique finie.
\end{defi}

On montre \cite{MP2} qu'il suffit de le v\'erifier pour une cohomologie $\ell$-adique si $k$ est de caract\'eristique nulle, et pour la cohomologie cristalline si $k$ est de caract\'eristique finie. De plus, pour une famille de tels endomorphismes, \^etre sp\'ecial est une propri\'et\'e ouverte et ferm\'ee.
On note $\LEnd(A)$ l'ensemble des endomorphismes sp\'eciaux de $A$.

La proposition suivante est due \`a Madapusi Pera \cite[Proposition 4.17 (4)]{MP} :

\begin{prop}\label{ident}
Si $A$ est la vari\'et\'e de Kuga-Satake associ\'ee \`a la surface K3 polaris\'ee $(X,\xi)$, il y a un isomorphisme $\Pic(X)^{\perp\xi}\xrightarrow{\sim}  \LEnd(A)$ compatible aux r\'ealisations cohomologiques. 
\end{prop}

\noindent{\sc Preuve} ---
Les arguments sont analogues \`a ceux des deux paragraphes pr\'ec\'edents. Expliquons-les rapidement, surtout pour mettre en \'evidence que le th\'eor\`eme de Torelli \ref{Torelli} est utilis\'e de mani\`ere cruciale.
Comme la preuve est plus facile en caract\'eristique nulle, supposons $X$ de caract\'eristique finie $p$.

Si $\zeta\in \Pic(X)^{\perp\xi}$, on peut relever $(X,\xi,\zeta)$ en $(\mathcal{X},\tilde{\xi},\tilde{\zeta})$ en caract\'eristique nulle \cite[Proposition A.1]{FM}, noter $\mathcal{A}$ la vari\'et\'e de Kuga-Satake associ\'ee et plonger son corps de d\'efinition dans $\mathbb{C}$. La classe $[\tilde{\zeta}_{\mathbb{C}}]_{\B}\in P^2_{\B}(\mathcal{X}_{\mathbb{C}})$ est de Hodge, induit une classe de Hodge dans $\End(H^1_{\B}(\mathcal{A}_{\mathbb{C}}))$, qui correspond \`a un endomorphisme sp\'ecial de $\mathcal{A}_{\mathbb{C}}$, qu'on peut sp\'ecialiser en un endomorphisme sp\'ecial de $A$.

R\'eciproquement, soit $f\in \LEnd(A)$ un endomorphisme sp\'ecial. L'\'etude des d\'eformations des endomorphismes sp\'eciaux, effectu\'ee \`a l'aide des th\'eories de Serre-Tate et de Grothendieck-Messing, montre qu'il est possible de relever le point $[A]\in \mathpzc{Sh}'_d(k)$ en $\mathcal{A}$ en caract\'eristique nulle de sorte que l'endomorphisme sp\'ecial $f$ se rel\`eve en un endomorphisme sp\'ecial $\tilde{f}$.
Le th\'eor\`eme \ref{Torelli} montre que ce rel\`evement induit un rel\`evement compatible $\mathcal{X}$ de la surface K3 $X$. La classe $[\tilde{f}_\mathbb{C}]_{\B}$ est de Hodge. Comme $\tilde{f}$ est sp\'ecial, elle correspond \`a une classe de Hodge dans $P^2_{\B}(\mathcal{X}_{\mathbb{C}})$. C'est la classe d'un fibr\'e en droites sur $\mathcal{X}_{\mathbb{C}}$, qu'on peut sp\'ecialiser sur $X$.
 \qed

\begin{rema}\label{quadraspecial}
La construction de Kuga-Satake sur $\mathbb{C}$ expliqu\'ee au paragraphe \ref{KSnaive} montre que si $\zeta\in \Pic(X)^{\perp\xi}$, l'endomorphisme sp\'ecial $f$ associ\'e satisfait $f\circ f=(\zeta^2)\Id$.
\end{rema}

On en d\'eduit imm\'ediatement :

\begin{coro}\label{corospecial}
Pour prouver le th\'eor\`eme \ref{thTate2} en caract\'eristique finie, il suffit de montrer que $\LEnd(A)\otimes\mathbb{Q}_{\ell}\to L_{\ell}(A)^{\Gamma_k}$ est un isomorphisme
\end{coro}

Montrons, en suivant \cite[Corollary 5.11]{MP}, que nous pouvons maintenant nous ramener au cas des corps finis, \`a l'aide de  la conjecture de Tate pour les vari\'et\'es ab\'eliennes sur les corps de type fini de caract\'eristique finie, prouv\'ee par Zarhin.

\begin{prop}\label{tf}
Pour prouver le th\'eor\`eme \ref{thTate2} en caract\'eristique finie, il suffit de traiter le cas o\`u $k$ est fini.
\end{prop}

\noindent{\sc Preuve} ---
\'Ecrivons $k$ comme corps des fonctions d'une vari\'et\'e int\`egre $B$ sur un corps fini. Quitte \`a r\'etr\'ecir $B$, on peut supposer que $X$ et $A$ ont des mod\`eles propres et lisses $\mathcal{X}$ et $\mathcal{A}$ sur $B$. Choisissons $b\in B(\mathbb{F}_q)$ un point ferm\'e de $B$. 

 Soient $\alpha\in L_\ell(A)^{\Gamma_k}$ une classe de cohomologie et $\alpha_b\in L_\ell(\mathcal{A}_b)^{\Gamma_{\mathbb{F}_q}}$ sa sp\'ecialisation. D'une part, par le th\'eor\`eme de Zarhin \cite[Corollary 2]{Zarhin2}, $\alpha$
est $\mathbb{Q}_{\ell}$-combinaison lin\'eaire de classes d'endomorphismes de $A$. D'autre part, par hypoth\`ese, 
 $\alpha_b$ est $\mathbb{Q}_{\ell}$-combinaison lin\'eaire de classes d'endomorphismes sp\'eciaux de $\mathcal{A}_b$. Par injectivit\'e des applications classe de cycle et 
 sp\'ecialisation, $\alpha$ est $\mathbb{Q}_{\ell}$-combinaison lin\'eaire de classes d'endomorphismes de $A$ qui sont sp\'eciaux sur $\mathcal{A}_b$, donc d'endomorphismes sp\'eciaux de $A$. 
 \qed

\subsection{R\'eduction modulo $p$ de vari\'et\'es de Shimura}\label{TMP}

Dans ce paragraphe, nous finissons de pr\'esenter la preuve par Madapusi Pera du th\'eor\`eme \ref{thTate2} en caract\'eristique finie. Par la proposition \ref{tf}, on suppose que $k$ est fini.

 Par le corollaire \ref{corospecial}, il suffit de montrer que toute classe dans $L_{\ell}(A)^{\Gamma_k}$ est \mbox{$\mathbb{Q}_\ell$-combinaison} lin\'eaire de classes d'endomorphismes sp\'eciaux.
Malheureusement, la conjecture de Tate pour les vari\'et\'es ab\'eliennes, qui montre seulement qu'elle est $\mathbb{Q}_\ell$-combinaison lin\'eaire de classes d'endomorphismes, ne permet pas de conclure.

\medskip

Cependant, cet \'enonc\'e ne fait plus intervenir de surfaces K3, mais seulement des vari\'et\'es ab\'eliennes obtenues par la construction de Kuga-Satake et leurs endomorphismes sp\'eciaux : c'est un probl\`eme interne \`a la th\'eorie des vari\'et\'es de Shimura de type $\GSpin$. La preuve de Madapusi Pera tire parti de ce cadre plus g\'en\'eral. 

Pour cela, plongeons $L_d\otimes\mathbb{Z}_{(p)}$ comme un sous-r\'eseau primitif propre d'un $\mathbb{Z}_{(p)}$-r\'eseau $M$ de signature $(2,m)$ et de discriminant premier \`a $p$ \cite[Lemma 6.8]{MP2}. On obtient un morphisme de vari\'et\'es de Shimura $\Sh'_d\to\Sh'_M$, qui induit un morphisme $f:\mathpzc{Sh}'_d\to\mathpzc{Sh}'_{M}$ entre leurs mod\`eles canoniques lisses en $p$. On remarque alors \cite[Lemma 5.6]{MP} que l'\'enonc\'e qu'on cherche \`a montrer pour $[A]\in\mathpzc{Sh}'_d(\ovFp)$ est une cons\'equence formelle de l'\'enonc\'e analogue pour $f([A])$. De plus, cet \'enonc\'e analogue est plus simple pour deux raisons : d'une part, le discriminant est maintenant premier \`a $p$, et d'autre part, l'orthogonal de $L_d$ dans $M$ se plonge dans $\LEnd(f([A]))$, qui est donc non trivial.

Pour ne pas alourdir les notations, nous supposons dans la suite que ces conditions \'etaient d\'ej\`a v\'erifi\'ees sans qu'il ait \'et\'e n\'ecessaire d'introduire le sur-r\'eseau $M$ : on a donc $p\nmid d$, et $\LEnd(A)\neq 0$.

\medskip

   Il est remarquable qu'il ait \'et\'e possible de forcer l'existence d'un endomorphisme sp\'ecial de mani\`ere aussi formelle. En effet, tous les autres seront construits \`a partir de celui-ci. Pour cela, Madapusi Pera applique des r\'esultats r\'ecents de Kisin \cite{Kisinnew} sur les r\'eductions modulo $p$ de vari\'et\'es de Shimura de type Hodge. Pour les d\'ecrire, nous nous pla\c{c}ons bri\`evement dans ce cadre plus g\'en\'eral.
  
  Soit $\Sh_K(G,X)$ une vari\'et\'e de Shimura de corps reflex $E$. On la suppose de type Hodge : il existe une immersion de vari\'et\'es de Shimura $\Sh_K(G,X)\to \Ab$ de but un champ de modules de vari\'et\'es ab\'eliennes polaris\'ees. Si $\mathfrak{p}$ est une place de $E$ au-dessus de $p$ qui est hypersp\'eciale\footnote{Dans notre application, c'est la condition $p\nmid d$.}, Kisin \cite[Theorem 2.3.8]{Kisin} construit un mod\`ele entier $\mathpzc{Sh}_K(G,X)$ de $\Sh_K(G,X)$ canonique lisse en $\mathfrak{p}$, \`a l'aide de la normalisation de son adh\'erence dans le mod\`ele entier naturel (c'est-\`a-dire modulaire) $\mathpzc{Ab}$ de $\Ab$ (dans notre cas, c'est ce que nous avons d\'ej\`a expliqu\'e au paragraphe \ref{parentier}). 
  
   Kisin donne une interpr\'etation modulaire, en un sens faible, de ce mod\`ele entier \cite[Proposition 1.3.9]{Kisinnew} : de la m\^eme mani\`ere que  $\Sh_K(G,X)$ s'interpr\`ete comme un espace de modules de vari\'et\'es ab\'eliennes polaris\'ees dont certaines classes de cohomologie sont de Hodge, les $\ovFp$-points de $\mathpzc{Sh}_K(G,X)$ correspondent \`a des vari\'et\'es ab\'eliennes sur $\ovFp$ munies de certaines classes de cohomologie distingu\'ees (en cohomologie $\ell$-adique et cristalline). L'ensemble $\mathpzc{Sh}_K(G,X)(\ovFp)$ est alors naturellement partitionn\'e en classes d'isog\'enie (o\`u l'on requiert que les isog\'enies pr\'eservent les classes de cohomologie distingu\'ees) \cite[Proposition 1.4.15]{Kisinnew}.
   
   Soient $x\in\mathpzc{Sh}_K(G,X)(\ovFp) $ et $A$ la vari\'et\'e ab\'elienne associ\'ee. Si $\ell\neq p$, le sous-groupe de $\GL(H^1_{\et}(A,\mathbb{Q}_{\ell}))$ pr\'eservant les classes de cohomologie distingu\'ees est isomorphe \`a~$G_{\mathbb{Q}_{\ell}}$. On note $I_{\ell}\subset G_{\mathbb{Q}_{\ell}}$ le centralisateur d'une puissance assez grande du Frobenius $\ell$-adique. Proc\'edant de m\^eme avec la cohomologie cristalline, on d\'efinit un sous-groupe $I_p\subset G_{\Frac(W(\ovFp))}$. 
Notons enfin $I$ le $\mathbb{Q}$-groupe alg\'ebrique des $\mathbb{Q}$-automorphismes de $A$ qui agissent \`a travers $I_{\ell}$ (resp. $I_p$) en cohomologie $\ell$-adique (resp. cristalline).
   
   Le th\'eor\`eme de Kisin que nous allons appliquer est \cite[Corollary 2.1.7]{Kisinnew} :
   
\begin{theo}\label{thkisin}
Il existe $\ell\neq p$ tel que $I_{\mathbb{Q}_{\ell}}\to I_{\ell}$ soit un isomorphisme\footnote{Kisin \cite[Corollary 2.3.2]{Kisinnew} en d\'eduit que c'est vrai pour tout $\ell$, mais nous n'en aurons pas besoin.}.
 \end{theo}

\noindent{\sc Id\'ee de la preuve} ---
Le groupe $I$ pr\'eserve \`a scalaire pr\`es la polarisation de $A$. Par positivit\'e de l'involution de Rosati associ\'ee, $I_{\mathbb{R}}$ est extension d'un groupe compact par $\mathbb{G}_{m,\mathbb{R}}$.
On en d\'eduit que $I$ est r\'eductif, donc que le quotient $Q:=I_{\ell}/I_{\mathbb{Q}_{\ell}}$ est affine. Comme le Frobenius agit de mani\`ere semi-simple sur $H^1_{\et}(A,\mathbb{Q}_{\ell})$, $I_{\ell}$ est le centralisateur d'un tore dans $G_{\mathbb{Q}_{\ell}}$, donc un sous-groupe de Levi de $G_{\mathbb{Q}_{\ell}}$. En particulier, il est g\'eom\'etriquement connexe et $Q$ est g\'eom\'etriquement connexe.

 Pour conclure, il reste \`a montrer que $Q$ est propre. Par un lemme de th\'eorie des groupes alg\'ebriques \cite[Lemma 2.1.8]{Kisinnew}, sous l'hypoth\`ese que $I_{\ell}$ est d\'eploy\'e, il suffit de montrer que $I_{\ell}(\mathbb{Q}_{\ell})/I(\mathbb{Q}_{\ell})$ est compact pour la topologie $\ell$-adique.
 \`A l'aide du th\'eor\`eme de Chebotarev, on assure cette hypoth\`ese en choisissant $\ell$ tel que $G_{\mathbb{Q}_{\ell}}$ soit d\'eploy\'e  et
tel que les valeurs propres du Frobenius agissant sur $H^1_{\et}(A,\mathbb{Q}_{\ell})$ appartiennent \`a $\mathbb{Q}_{\ell}$.

La description que donne Kisin des classes d'isog\'enie montre l'existence, pour un $q$ bien choisi, d'une application $I(\mathbb{Q})\backslash I_{\ell}(\mathbb{Q}_{\ell})/(I_{\ell}(\mathbb{Q}_{\ell})\cap K_{\ell})\to \mathpzc{Sh}_K(G,X)(\mathbb{F}_q)$. 
L'inter\-pr\'etation modulaire de $\mathpzc{Sh}_K(G,X)(\ovFp)$ lui permet de montrer (essentiellement dans la preuve de \cite[Proposition 2.1.3]{Kisinnew}) l'injectivit\'e de cette application. C'est l'\'etape cruciale : on construit des \'el\'ements de $I(\mathbb{Q})$, donc des automorphismes de $A$.

 Il s'ensuit que $I(\mathbb{Q_{\ell}})\backslash I_{\ell}(\mathbb{Q}_{\ell})/(I_{\ell}(\mathbb{Q}_{\ell})\cap K_{\ell})$ est fini : $I_{\ell}(\mathbb{Q}_{\ell})/I(\mathbb{Q}_{\ell})$ est donc compact.
 \qed

\medskip
   
   On peut maintenant conclure en suivant la preuve de \cite[Theorem 5.4]{MP}. Commen\c{c}ons par choisir $\ell$ comme dans le th\'eor\`eme \ref{thkisin}. Dans le cas particulier des vari\'et\'es de Shimura de type $\GSpin$ que nous consid\'erions, la d\'efinition de $I_{\ell}$ montre que, quitte \`a remplacer $k$ par une extension finie, $I_{\ell}$ agit sur $H_{\ell}^{\otimes (1,1)}(A)$ en pr\'eservant $L_{\ell}(A)^{\Gamma_k}$, et que son action sur $L_{\ell}(A)^{\Gamma_k}$ est irr\'eductible \cite[Lemma 5.8]{MP}. Par le th\'eor\`eme \ref{thkisin}, l'action de $I$ sur $L_{\ell}(A)^{\Gamma_k}$ est irr\'eductible. Comme $I$ agit par automorphismes de $A$, cette action pr\'eserve l'espace $\LEnd(A)\otimes\mathbb{Q}_{\ell}$ des endomorphismes sp\'eciaux. 
 Nous nous \'etions ramen\'es au cas o\`u cet espace est non nul, de sorte que l'irr\'eductibilit\'e de l'action montre que $\LEnd(A)\otimes\mathbb{Q}_{\ell}\to L_{\ell}(A)^{\Gamma_k}$ est un isomorphisme.
 
 Nous avons montr\'e la conjecture de Tate pour $X$ et un $\ell$ particulier. On peut la reformuler comme suit : le nombre de Picard de $X$ est \'egal \`a la multiplicit\'e de la valeur propre $1$ pour le Frobenius $\ell$-adique. Par ind\'ependance de $\ell$ du polyn\^ome caract\'eristique du Frobenius $\ell$-adique \cite{Weil1}, et semi-simplicit\'e de l'action de celui-ci \cite[Corollaire 1.10]{K3rel}, on en d\'eduit imm\'ediatement le m\^eme \'enonc\'e pour tout nombre premier $\ell\neq p$.
 
 \medskip

La parent\'e de l'argument ci-dessus avec la preuve par Tate \cite{Tateendo} de la conjecture de Tate pour les vari\'et\'es ab\'eliennes sur les corps finis est frappante : celle-ci repose \'egalement sur un r\'esultat de finitude de classes d'isog\'enie, et fait aussi, dans un premier temps, le choix d'un nombre premier $\ell$ particulier.

\subsection{Utilisation de formes automorphes}\label{MauCha}

Le but de ce paragraphe est de prouver le th\'eor\`eme suivant :

\begin{theo}\label{thTate3}
Soit $X$ une surface K3 supersinguli\`ere sur un corps alg\'ebriquement clos $k$ de carac\-t\'eristique~$p \neq 2$. Si $X$ poss\`ede une polarisation de degr\'e $2d$ premier \`a $p$, $\rho(X)=22$.
\end{theo}

Sous l'hypoth\`ese plus forte $p>2d+4$, ce th\'eor\`eme a \'et\'e obtenu par Maulik \cite{Maulik}. Le cas $p\nmid 2d$ est d\^u \`a Charles \cite{Charles}, qui reprend et am\'eliore la strat\'egie de Maulik.


L'id\'ee remonte \`a \cite[Theorem 4]{RSZ} : par un th\'eor\`eme d'Artin \cite[Theorem 1.1]{Artin}, il suffit de montrer le th\'eor\`eme \ref{thTate3} pour une surface K3 particuli\`ere dans chaque composante irr\'eductible $C\subset\widetilde{\mathfrak{M}}_{2d}'^{\lisse}$ du lieu supersingulier. Comme ce th\'eor\`eme est connu dans le cas elliptique \cite[Theorem 1.7]{Artin}, il suffit de montrer que $C$ intersecte le lieu elliptique\footnote{Quand $p=3$, on doit inclure dans le lieu elliptique les surfaces K3 munies de fibrations quasi-elliptiques. Celles-ci v\'erifient la conjecture de Tate car elles sont unirationnelles, donc Shioda-supersinguli\`eres \cite[Corollary 2]{Shiodass}.}. Pour cela, nous allons combiner deux ingr\'edients : une propri\'et\'e d'amplitude du lieu elliptique, et un r\'esultat de propret\'e du lieu supersingulier. 

\medskip

Le premier, d\^u \`a Maulik \cite[\S 3]{Maulik}, repose sur des constructions de formes automorphes, et exploite \`a cet effet des travaux de Borcherds \cite{Borcherds}.
Commen\c{c}ons par travailler sur $\mathbb{C}$, en conservant les notations du paragraphe \ref{parperiodes}. Pour tout $v\in L_{d}^\vee$, on dispose d'un diviseur $v^{\perp}\subset \Omega$. Si $n\in \mathbb{Q}^{>0}$ et $\gamma\in L_d^\vee/L_d$,  le diviseur $\sum_{\langle v,v\rangle=-2n, [v]=\gamma}v^{\perp}$ est $\Gamma$-invariant et descend en un diviseur $y_{n,\gamma}$ sur $\Omega/\Gamma=\Sh_{d,\mathbb{C}}$ : c'est un \textbf{diviseur de Heegner}.
Par convention, on d\'efinit $y_{0,\gamma}=0$ si $\gamma\neq 0$, et $[y_{0,0}]=\lambda_{\mathbb{C}}$, o\`u $\lambda_{\mathbb{C}}$ est le fibr\'e de Hodge sur $\Sh_{d,\mathbb{C}}$.

Si $k\geq 1$, le th\'eor\`eme des classes $(1,1)$ de Lefschetz montre que le groupe de Picard d'une surface K3 complexe dont les p\'eriodes appartiennent au diviseur $y_{dk^2,0}$ contient un sous-r\'eseau isomorphe \`a $\left(\begin{matrix} 2d & 0\\ 0  & -2dk^2\end{matrix}\right)$. Comme cette forme quadratique est isotrope, une telle surface K3 est elliptique.

Supposons maintenant par l'absurde que $\lambda_{\mathbb{C}}$ ne s'\'ecrit pas comme combinaison lin\'eaire dans $\Pic_{\mathbb{Q}}(\Omega/\Gamma)$ des $[y_{dk^2,0}]$. Alors il existe une forme lin\'eaire $\alpha$ sur  $\Pic_{\mathbb{Q}}(\Omega/\Gamma)$ telle que $\alpha([y_{dk^2,0}])=0$ pour tout $k$, mais $\alpha(\lambda)\neq 0$. Un th\'eor\`eme de Borcherds \cite[Theorem 4.5]{Borcherds} pr\'ecis\'e par McGraw \cite{McGraw} montre que la s\'erie formelle :
$$\Phi_{\alpha}(q)=\sum_{n\in \mathbb{Q}^{\geq 0}}\sum_{\gamma\in L_d^\vee/L_d}\alpha([y_{n,\gamma}])v_{\gamma}q^n\in \mathbb{Q}[L_d^\vee/L_d][[q^{\frac{1}{4d}}]]$$ est le d\'eveloppement de Fourier d'une forme modulaire vectorielle \`a valeurs dans $\mathbb{C}[L_d^\vee/L_d]$ et de poids demi-entier $21/2$.
Comme $\alpha([y_{0,0}])\neq 0$, on peut \'ecrire $\Phi_{\alpha}(q)$ comme somme d'une forme modulaire non nulle construite \`a l'aide d'une s\'erie d'Eisenstein, et d'une forme modulaire cuspidale. Le coefficient en $v_0q^{dk^2}$ de la premi\`ere cro\^it plus rapidement avec $k$ que celui de la seconde. Cela contredit, pour $k\gg 0$, le fait que $\alpha([y_{dk^2,0}])=0$.

Nous avons donc montr\'e \cite[Theorem 3.1]{Maulik} l'existence d'un diviseur de Cartier (non n\'ecessairement effectif) $D$ sur $\Sh_{d,\mathbb{C}}$ support\'e sur le lieu elliptique,
et dont le fibr\'e en droites associ\'e est une puissance du fibr\'e de Hodge $\lambda_{\mathbb{C}}$. 
Le diviseur $D$ (dont les composantes irr\'eductibles sont en fait des sous-vari\'et\'es de Shimura de $\Sh_{d,\mathbb{C}}$) est d\'efini sur un corps de nombres. Quitte \`a le remplacer par la somme de ses conjugu\'es sous Galois, on peut le supposer d\'efini sur $\mathbb{Q}$. On note encore $D$ son pull-back \`a $\Sh'_{d}$. Le r\'esultat que nous utiliserons est l'amplitude de sa r\'eduction modulo $p$ :

\begin{prop}
L'adh\'erence $\overline{D}$ de $D$ dans $\mathpzc{Sh}'_d$ est ample en restriction \`a $\mathpzc{Sh}'_{d,\mathbb{F}_p}$.
\end{prop}

\noindent{\sc Preuve} --- Notons $\mu$ le fibr\'e ample naturel sur $\mathpzc{Ab}$ : le d\'eterminant du fibr\'e de Hodge \cite[Theorem V.2.3]{FC}. Comme $\iota$ est quasi-fini, $\iota^*\mu$ est un fibr\'e ample sur $\mathpzc{Sh}'_d$. 
Maulik montre \cite[Proposition 5.8]{Maulik} que $\lambda_{\mathbb{C}}$ est proportionnel \`a $\iota_{\mathbb{C}}^*\mu_{\mathbb{C}}$. Des multiples convenables des fibr\'es en droites  $\mathcal{O}(\overline{D})$ et $\iota^*\mu$ sur $\mathpzc{Sh}'_{d,\mathbb{Z}_{(p)}}$ co\"incident donc en restriction \`a  $\mathpzc{Sh}'_{d,\mathbb{Q}}$ :
ils diff\`erent d'un diviseur de Cartier $\Delta$ support\'e sur $\mathpzc{Sh}'_{d,\mathbb{F}_p}$. Comme $\mathpzc{Sh}'_{d,\mathbb{Z}_{(p)}}$ est lisse, $\Delta$ est trivial en restriction \`a $\mathpzc{Sh}'_{d,\mathbb{F}_p}$ \cite[Lemma 5.12]{Maulik}, et $\mathcal{O}(\overline{D})|_{\mathpzc{Sh}'_{d,\mathbb{F}_p}}$ est bien ample.
\qed

\medskip

Le deuxi\`eme ingr\'edient est un r\'esultat de propret\'e du lieu supersingulier. Maulik suppose que $p>2d+4$. Il montre qu'une surface K3 polaris\'ee supersinguli\`ere sur le corps des fractions d'un anneau de valuation discr\`ete de caract\'eristique $p$, dont la polarisation est tr\`es ample,
a potentiellement bonne r\'eduction\footnote{Ce r\'esultat \'etait d\'ej\`a connu sous la conjecture de Tate : c'est le th\'eor\`eme de Rudakov et Shafarevich \cite{RS} que nous avons \'enonc\'e au th\'eor\`eme \ref{RudaSha}.} \cite[Theorem 4.1]{Maulik}. En effet, sous cette restriction sur la caract\'eristique, un th\'eor\`eme de Saito \cite{Saito} permet de construire un mod\`ele semistable. Maulik applique le programme du mod\`ele minimal en caract\'eristique $p$ de Kawamata \cite{Kawamata} \`a ce mod\`ele semistable, et v\'erifie que sous l'hypoth\`ese de supersingularit\'e, le mod\`ele obtenu par ce proc\'ed\'e est lisse.

Cela permet \`a Maulik de construire des familles compl\`etes non triviales de surfaces K3 supersinguli\`eres. Une difficult\'e est que ces familles pourraient ne pas \^etre polaris\'ees, mais Maulik est capable d'assurer la condition plus faible, et suffisante pour son argument, qu'elles soient quasi-polaris\'ees. Finalement, le th\'eor\`eme de Torelli \ref{Torelli} implique que cette famille n'est pas contract\'ee par l'application des p\'eriodes $P'$, et permet de montrer qu'elle rencontre le diviseur ample $\overline{D}|_{\mathpzc{Sh}'_{d,\mathbb{F}_p}}$, donc le lieu elliptique. Nous ne d\'evelopperons pas plus ces arguments, qui sont d\'etaill\'es dans \cite[\S 7]{Maulik}. 

\medskip

Quand $p\nmid 2d$, l'id\'ee de Charles est d'exploiter plut\^ot la propret\'e, prouv\'ee par Oort \cite[Theorem 1.1 (a)]{Oortcomplete}, du lieu supersingulier dans l'espace de modules de vari\'et\'es ab\'eliennes polaris\'ees $\mathpzc{Ab}_{\mathbb{F}_p}$\footnote{Oort montre un r\'esultat plus fort : le lieu des vari\'et\'es ab\'eliennes de $p$-rang nul est propre.}. Ce r\'esultat est plus facile, et valable en toute caract\'eristique : dans sa preuve, l'existence de mod\`eles de N\'eron sert de substitut \`a l'existence de mod\`eles semistables.
Qu'il soit utile dans notre situation est montr\'e par la proposition suivante \cite[Proposition 21]{Charles} :

\begin{prop}\label{ssss}
La vari\'et\'e de Kuga-Satake $A$ d'une surface K3 super\-singuli\`ere~$X$  est supersinguli\`ere.
\end{prop}

\noindent{\sc Preuve} ---
Nous devons introduire une variante de la construction de Kuga-Satake, qui est pr\'ecis\'ement celle consid\'er\'ee par Deligne \cite{WeilK3}. Reprenons temporairement les notations du paragraphe \ref{KSnaive}.  Si nous avions utilis\'e l'alg\`ebre de Clifford paire, nous aurions construit une vari\'et\'e ab\'elienne $A^+:=\Cl^+(P^2_{\B}(X,\mathbb{R}))/\Cl^+(P^2_{\B}(X,\mathbb{Z}))$, dite \textbf{vari\'et\'e de Kuga-Satake paire} de $X$, telle que $A$ est isog\`ene \`a $(A^{+})^{2}$. L'alg\`ebre de Clifford paire $\Cl^+(L_d)$ agit sur $A^+$ par multiplication \`a droite. On v\'erifie qu'on a alors un isomorphisme de structures de Hodge $\Cl^+(P^2_{\B}(X)(1))\simeq\End_{\Cl^+(L_d)}(H^1_{\B}(A^+))$ donn\'e par la multiplication \`a gauche. Des arguments analogues \`a ceux d\'ecrits au paragraphe \ref{parQ} permettent de d\'efinir $A^+$ et l'action de $\Cl^+(L_d)$ sur tout corps, et d'y \'etendre l'isomorphisme ci-dessus en cohomologie $\ell$-adique \cite[6.5, 6.6]{WeilK3}.


Prouvons maintenant la proposition.
Par sp\'ecialisation, on voit qu'il suffit de traiter le cas o\`u $X$ est d\'efini sur un corps fini $k$.  Comme $A$ est isog\`ene \`a $(A^{+})^{2}$, il suffit de montrer que $A^+$ est supersinguli\`ere.
 On a un isomorphisme de $\Gamma_k$-modules $\Cl^+(P^2_{\et}(X_{\bar{k}},\mathbb{Q}_{\ell}(1)))\simeq\End_{\Cl^+(L_d)}(H^1_{\et}(A^+_{\bar{k}},\mathbb{Q}_{\ell}))$.  Par la proposition \ref{defss} (iii), on peut supposer que le Frobenius agit trivialement sur $P^2_{\et}(X_{\bar{k}},\mathbb{Q}_{\ell}(1))$, donc sur $\End_{\Cl^+(L_d)}(H^1_{\et}(A^+_{\bar{k}},\mathbb{Q}_{\ell}))$. Cela signifie qu'il agit sur $H^1_{\et}(A^+_{\bar{k}},\mathbb{Q}_{\ell})$ par un \'el\'ement du bicommutant de $\Cl^+(L_d)$. Comparant avec la cohomologie de Betti, on voit que ce bicommutant est \'egal \`a $\Cl^+(L_d)\otimes \mathbb{Q}_{\ell}$. Comme le Frobenius commute \`a l'action de $\Cl^+(L_d)$, il agit via un \'el\'ement du centre de $\Cl^+(L_d)\otimes \mathbb{Q}_{\ell}$. Comme $\Cl^+(L_d)\otimes \mathbb{Q}_{\ell}$ est centrale simple\footnote{C'est la raison pour laquelle nous avons d\^u utiliser l'alg\`ebre de Clifford paire.} \cite[\S 9.4, Th\'eor\`eme~3]{Bourbaki}, il agit par homoth\'etie, et $A^+$ est donc supersinguli\`ere.
\qed
\medskip

 On peut maintenant achever la d\'emonstration du th\'eor\`eme \ref{thTate3} en suivant la preuve de \cite[Proposition 22]{Charles}. Rappelons que, par le th\'eor\`eme \ref{Torelli}, l'application des p\'eriodes $P'$ est une immersion ouverte. Fixons une composante irr\'eductible $C$ du lieu supersingulier, notons $\overline{C}$ l'adh\'erence de $C$ dans $\mathpzc{Sh}'_d$, $\bar{\eta}$ un point g\'en\'erique g\'eom\'etrique de $\overline{C}$ et $\mathcal{A}$ le sch\'ema ab\'elien sur $\overline{C}$ induit par $\iota$.

La construction de Kisin expliqu\'ee au paragraphe \ref{parentier} montre que $\iota$ est propre\footnote{C'est ici que nous utilisons que la polarisation est premi\`ere \`a la caract\'eristique.}. Comme le lieu supersingulier dans $\mathpzc{Ab}_{\mathbb{F}_p}$ est propre, son image r\'eciproque par $\iota$ l'est \'egalement. La proposition \ref{ssss} montre alors que $\overline{C}$ est propre.
Comme $C$ est de dimension~$>0$ \cite[Theorem 15]{O3}, $\overline{C}$ intersecte le diviseur ample $\overline{D}|_{\mathpzc{Sh}'_{d,\mathbb{F}_p}}$. Soit $x$ un $\ovFp$-point dans l'intersection. Il est possible que $x$ n'appartienne par \`a $C$, donc que $\mathcal{A}_x$ ne soit pas la vari\'et\'e de Kuga-Satake associ\'ee \`a une surface K3. Cependant, le choix de $D$ montre que $\mathcal{A}_x$ est sp\'ecialisation de la vari\'et\'e de Kuga-Satake d'une surface K3 particuli\`ere, elliptique, en caract\'eristique nulle. On en d\'eduit par la remarque \ref{quadraspecial} qu'il existe un entier $k\geq 1$ tel que $\mathcal{A}_x$ poss\`ede un endomorphisme sp\'ecial $\phi$ tel que $\phi\circ\phi=-2dk^2\Id$.

Comme $\mathcal{A}_{\bar{\eta}}$ est supersinguli\`ere, le rang de $\End(\mathcal{A}_{\bar{\eta}})$ est maximal, de sorte que le conoyau de l'application de sp\'ecialisation $\End(\mathcal{A}_{\bar{\eta}})\to \End(\mathcal{A}_x)$ est de torsion. Un multiple de $\phi$ se rel\`eve donc en un endomorphisme sp\'ecial $\psi\in\End(\mathcal{A}_{\bar{\eta}})$, qui satisfait $\psi\circ\psi=-2dk'^2\Id$ pour un certain entier $k'\geq 1$. La surface K3 associ\'ee \`a $\bar{\eta}$ poss\`ede, par la proposition \ref{ident} et la remarque \ref{quadraspecial}, un fibr\'e en droites $\zeta$ orthogonal \`a la polarisation tel que $\zeta^2=-2dk'^2$. Comme $\left(\begin{matrix} 2d & 0\\ 0  & -2dk'^2\end{matrix}\right)$ est isotrope, elle est elliptique. Cela conclut.

\section{Construction de courbes rationnelles}

Dans cette partie, nous expliquons la preuve du th\'eor\`eme suivant \cite{BHT, LL}. On trouvera \'egalement une exposition de ces r\'esultats dans \cite[\S 13.3]{HuyK3}.

\begin{theo}  
\label{BHTLL2}
Une surface K3 $X$ sur un corps alg\'ebriquement clos $k$ contient une infinit\'e de courbes rationnelles\footnote{Dans tout ce texte, une courbe rationnelle est une sous-vari\'et\'e int\`egre dont la normalisation est isomorphe \`a $\mathbb{P}^1$.} dans les deux cas suivants :
\begin{enumerate}[(i)]
\item Si $k$ est de caract\'eristique nulle et $X$ n'est pas d\'efinie sur $\ovQ$.
\item Si $k$ est de caract\'eristique $\neq 2$ et $\rho(X)$ est impair.
\end{enumerate}
\end{theo}

  La strat\'egie de la preuve est due \`a Bogomolov, Hassett et Tschinkel : ils font fonctionner celle-ci sous l'hypoth\`ese suppl\'ementaire que $X$ admet une polarisation de degr\'e~$2$. Li et Liedtke ont obtenu le cas g\'en\'eral en reprenant et en am\'eliorant leurs id\'ees.

  Les techniques de la preuve s'inscrivent dans la lign\'ee de travaux ant\'erieurs. Un th\'eor\`eme attribu\'e par Mori et Mukai \cite{MM1} \`a Bogomolov et Mumford montre l'existence d'une courbe rationnelle sur toute surface K3. Chen \cite{Chen}
a am\'elior\'e ce r\'esultat en prouvant qu'une surface K3 complexe tr\`es g\'en\'erale contient une infinit\'e de courbes rationnelles nodales. Ces r\'esultats, comme le th\'eor\`eme qui nous int\'eresse, reposent sur des arguments de d\'eformation. Ils exploitent la remarque suivante : on peut d\'eformer une courbe rationnelle $C$ sur une surface K3 \`a des d\'eformations de celle-ci d\`es que le fibr\'e en droites $\mathcal{O}(C)$ s'y d\'eforme. Ce principe sera notre principal outil de construction de courbes rationnelles.

  Les paragraphes \ref{parr1} et \ref{parr2} sont consacr\'es \`a ces techniques de d\'eformation. On d\'emontre (proposition \ref{defo}) un \'enonc\'e qui rend rigoureux le principe formul\'e ci-dessus. Un des apports de \cite{BHT} est qu'il est utile de consid\'erer des d\'eformations d'applications stables plut\^ot que de simples courbes rationnelles : c'est dans ce cadre que se place la proposition \ref{defo}. 

L'am\'elioration principale qui a permis \`a Li et Liedtke \cite{LL} d'aller au-del\`a de l'article \cite{BHT} r\'eside dans un choix astucieux des applications stables \`a d\'eformer, et est expliqu\'ee au paragraphe \ref{parr2} (proposition \ref{rigide}).

  Les paragraphes \ref{parr3} et \ref{parr4} pr\'esentent la preuve du th\'eor\`eme \ref{BHTLL2}. Le cas (i) (resp. (ii)) est un argument de sp\'ecialisation et d\'eformation en \'egale caract\'eristique nulle (resp. en caract\'eristique mixte). Les articles \cite{BHT} et \cite{LL} se placent dans le second cas, le premier en est une variante plus \'el\'ementaire, qui nous a \'et\'e signal\'ee par Charles. Dans les deux cas, il est ais\'e de construire des courbes rationnelles sur $X$, et la difficult\'e est de s'assurer qu'elles sont distinctes. Pour cela, on exploite le fait que $X$ poss\`ede une infinit\'e de sp\'ecialisations en caract\'eristique nulle (resp. en caract\'eristique finie) en lesquelles le nombre de Picard g\'eom\'etrique cro\^it. Cela r\'esulte d'un argument de th\'eorie de Hodge (resp. de la conjecture de Tate pour les surfaces K3 sur les corps finis). Si l'utilisation de la caract\'eristique finie
est la principale innovation de \cite{BHT}, les d\'etails de la preuve pr\'esent\'ee ici sont issus de \cite{LL}.
  
  Finalement, on discute au paragraphe \ref{parr5} les limites de cette strat\'egie.

\subsection{D\'eformation de courbes rationnelles}\label{parr1}

Commen\c{c}ons par \'enoncer le th\'eor\`eme de Bogomolov et Mumford (\cite[Appendix]{MM1}, \cite[Proposition 2.5]{BT}, \cite[Theorem 2.1]{LL}) mentionn\'e dans l'introduction, qui nous sera utile plus tard :

\begin{theo}\label{unecourbe}
Soient $X$ une surface K3 sur un corps alg\'ebriquement clos $k$ et $L$~un fibr\'e en droites effectif
sur $X$. Alors il existe un diviseur $D\in|L|$ dont toutes les composantes irr\'eductibles sont des courbes rationnelles.
\end{theo}

Ce th\'eor\`eme est tout d'abord d\'emontr\'e pour une surface K3 particuli\`ere (une surface de Kummer) en caract\'eristique nulle, en construisant le diviseur $D$ explicitement. On traite le cas d'une surface g\'en\'erale en caract\'eristique nulle en d\'eformant ce diviseur. Un argument de sp\'ecialisation montre alors le th\'eor\`eme pour toute surface K3 en toute caract\'eristique. Comme le diviseur explicite $D$ utilis\'e dans \cite{MM1} est nodal, \cite{HuyK3} remarque qu'on obtient en fait  :

\begin{theo}\label{nodale}
Soit  $(X,L)$ une surface K3 polaris\'ee g\'en\'erale sur un corps alg\'ebri\-que\-ment clos de caract\'eristique nulle. Alors $X$ contient une courbe rationnelle nodale.
\end{theo}

Ces arguments de d\'eformation seront importants pour la preuve du th\'eor\`eme \ref{BHTLL2}. Donnons un exemple d'\'enonc\'e, en suivant \cite[Theorem 18]{BHT} (voir aussi \cite[\S 2]{HuyKe}). On se place dans le contexte des applications stables de Kontsevich \cite{FP} : $\overline{\mathcal{M}}_0(X)$ d\'esigne le champ de modules des applications stables de genre $0$ sur $X$.

\begin{prop}\label{defo}
Soient $S$ un sch\'ema int\`egre de type fini sur un corps $k$ alg\'ebrique\-ment clos de caract\'eristique nulle, $(X,L)\to S$ une famille de surfaces $K3$ polaris\'ees, et $s\in S$. Soit $[f:C\to X_s]\in \overline{\mathcal{M}} :=\overline{\mathcal M}_0(X/S)$ une application stable de genre $0$ telle que $f_*[C]\in |L_s|$. On suppose que $f$ est \textbf{rigide} au sens o\`u la composante connexe de la fibre de $\overline{\mathcal{M}} \to S$ contenant $[f]$ est de dimension nulle. Alors il existe une composante de $\overline{\mathcal{M}}$ contenant $[f]$ qui domine $S$.  
\end{prop}

\noindent{\sc Preuve} --- 
   Soit $\mathcal{X}\to \Delta:=\Spf k[[x_1,\dots,x_{20}]]$ la d\'eformation verselle de $X_{s}$.  Par \cite[Lemma 11]{BHT}, la th\'eorie des d\'eformations de $f$ est contr\^ol\'ee par un faisceau coh\'erent $\mathcal{N}_f$\footnote{Quand $C$ est lisse et $f$ immersive, il s'agit du fibr\'e normal usuel.}
sur $C$ : les d\'eformations infinit\'esimales et les obstructions sont respectivement donn\'ees par $H^0(C,\mathcal{N}_f)$ et $H^1(C,\mathcal{N}_f)$. De plus, la description qui y est donn\'ee de $\mathcal{N}_f$, notamment la suite exacte en bas de la page 541, permet de calculer $\chi(C,\mathcal{N}_f)=-1$.
  
    La th\'eorie des d\'eformations montre alors que $\overline{\mathcal M}_0(\mathcal X/\Delta)$ est de dimension relative $\geq \chi(C, \mathcal{N}_f)=-1$ sur $\Delta$ au voisinage de $[f]$ (voir par exemple l'argument de \cite[I.2.15.1]{Kollar} qui, \'ecrit pour le sch\'ema de Hilbert, est g\'en\'eral et s'applique). Par l'hypoth\`ese de rigidit\'e, l'image des composantes irr\'eductibles de $\overline{\mathcal M}_0(\mathcal X/\Delta)$ passant par $[f]$ est de codimension $\leq 1$ dans $\Delta$. Comme le lieu o\`u le fibr\'e en droites $L_{s}$ se d\'eforme est un diviseur irr\'eductible $H\subset \Delta$
, et que cette image est incluse dans $H$, elle est n\'ecessairement \'egale \`a $H$.
   
   Comme l'image de $\Spf(\widehat{\mathcal{O}}_{S,s})\to \Delta$ est incluse dans $H$, on en d\'eduit l'existence d'une composante irr\'eductible de $\overline{\mathcal{M}}$ contenant $[f]$ dont l'image contient un voisinage de $s$, et qui domine donc $S$.
\qed

\begin{rema}
Les travaux de Bloch \cite{Bloch}, Voisin \cite{Voreg} et Ran \cite{Ran} sur l'application de semi-r\'egularit\'e permettent de comprendre la proposition \ref{defo} du point de vue de la th\'eorie de Hodge.
\end{rema}

\begin{rema}\label{deforem}
Quand $S$ est de caract\'eristique finie ou mixte, on peut obtenir des analogues de la proposition \ref{defo}, en consid\'erant la d\'eformation verselle de $X_s$ au-dessus de $\Spf(W(\kappa(s))[[x_1,\dots,x_{20}]])$ \cite[Theorem 18]{BHT}. Il convient de prendre garde que le dernier argument ne fonctionne pas en g\'en\'eral, car il est possible que $H$ ne soit pas irr\'eductible.
Cette difficult\'e dispara\^it si $X_s$ est ordinaire ou si $L_s$ est une polarisation primitive. 
\end{rema}

\subsection{Applications stables rigides}\label{parr2}

  La difficult\'e qu'ont rencontr\'ee Bogomolov, Hassett et Tschinkel pour prouver le th\'eor\`eme \ref{BHTLL2} sans hypoth\`ese sur la polarisation de $X$ est la suivante : ils souhaitaient d\'eformer une r\'eunion connexe $D$ de courbes rationnelles sur une surface K3 $X$ \`a des d\'eformations de $X$ en utilisant la proposition \ref{defo}. Si toutes les courbes rationnelles dans~$D$ sont distinctes, ce n'est pas difficile car toute courbe stable $f$ de genre $0$ ayant $D$ pour cycle associ\'e est rigide si $X$ n'est pas unir\'egl\'ee. En revanche, s'il y a des multiplicit\'es dans $D$, ce n'est pas automatique, m\^eme si $X$ n'est pas unir\'egl\'ee : $f$ pourrait avoir des d\'eformations non triviales \`a cycle constant.
  
  Li et Liedtke ont propos\'e la solution suivante \`a cette difficult\'e technique \cite[Theorem 3.9]{LL} :

\begin{prop}\label{rigide}
Soit $X$ une surface K3 non unir\'egl\'ee sur un corps alg\'ebri\-que\-ment clos. Soient $D_1,\dots,D_n, R$ des courbes rationnelles sur $X$. On suppose que $R$ est ample et nodale.
Alors il existe $k\geq 0$ et une courbe stable $[f]\in\overline{\mathcal{M}}_0(X)$ g\'en\'eriquement non ramifi\'ee et rigide dont le cycle associ\'e est $D_1+\dots+D_n+kR$.
\end{prop}

La preuve de cette proposition est combinatoire et astucieuse. Il faut recoller les normalisations des courbes $D_1,\dots,D_n, R$ pour obtenir une courbe stable de genre $0$, et le faire soigneusement, de sorte que celle-ci soit rigide. Nous ne recopierons pas ici les arguments de \cite[Lemma 3.6, Theorem 3.9]{LL}.

Bien s\^ur, une difficult\'e pour appliquer la proposition \ref{rigide} est de produire une courbe rationnelle nodale $R$. On dispose heureusement du th\'eor\`eme \ref{nodale} !

\subsection{Surfaces K3 qui ne sont pas d\'efinies sur $\ovQ$}\label{parr3}

Prouvons maintenant le th\'eor\`eme \ref{BHTLL2} (i).  Dans ce cas, la preuve est purement en caract\'eristique nulle : c'est l\'eg\`ement plus simple que le cadre de caract\'eristique mixte dans lequel se placent \cite{BHT} et \cite{LL}, et qui sera expliqu\'e au paragraphe suivant.

\begin{theo}\label{th1}
Soit $X$ une surface K3 sur un corps alg\'ebriquement clos de carac\-t\'eristique nulle qui n'est pas d\'efinie sur $\ovQ$. Alors $X$ contient une infinit\'e de courbes rationnelles.
\end{theo}

\noindent{\sc Preuve} --- 
Soit $A$ une polarisation primitive de degr\'e $2d$ sur $X$. Soient $U\to \mathfrak{M}_{2d,\ovQ}$ une carte \'etale int\`egre\footnote{Comme dans \cite{LL}, introduire $U$ permet d'\'eviter de manipuler des champs de modules d'applications stables sur une base champ\^etre.} du champ de modules des surfaces K3 munies d'une polarisation primitive de degr\'e $2d$, $(\mathcal{X},\mathcal{A})$ la famille universelle sur $U$, $\bar{\eta}$ un point g\'eom\'etrique de $U$ tel que $(X,A)\simeq(\mathcal{X}_{\bar{\eta}},\mathcal{A}_{\bar{\eta}})$, $S$ l'adh\'erence de $\bar{\eta}$ dans $U$ et $\bar{\xi}$ un point g\'en\'erique g\'eom\'etrique de $U$.

Comme $X$ n'est pas d\'efinie sur $\ovQ$, $\dim(S)>0$. Un th\'eor\`eme de Green \cite[Proposition 17.20]{Voisin}, appliqu\'e \`a la variation de structure de Hodge sur $S_{\mathbb{C}}$ donn\'ee par la cohomologie transcendante (i.e. l'orthogonal de $\Pic(X)$) de la famille de surfaces K3 $\mathcal{X}_{S_{\mathbb{C}}}\to S_{\mathbb{C}}$, montre que l'ensemble $\Sigma$ des points ferm\'es $s\in S$ tels que $\rho(\mathcal X_s)>\rho(X)$ est dense dans $S$. Par le th\'eor\`eme \ref{unecourbe}, pour tout $s\in \Sigma$, il existe une courbe rationnelle $C_s\subset \mathcal X_s$ dont la classe n'appartient pas \`a $\Pic(X)$. On note $f_s:\mathbb{P}^1\to\mathcal{X}_s$ sa normalisation.

Fixons un entier $N>0$. Si $\Sigma':=\{s\in \Sigma|\mathcal{A}_s.C_s>N\}$ n'\'etait pas dense dans $S$, les $([f_s])_{s\in\Sigma\setminus\Sigma'}$ appartiendraient tous au $S$-sch\'ema de type fini $M:=\Mor_{S,<N}(\mathbb{P}^1,\mathcal{X})$ param\'etrant les morphismes de degr\'e $<N$ de $\mathbb{P}^1$ vers $\mathcal{X}$. Une des composantes irr\'eductibles de $M$ contenant un de ces $[f_s]$ dominerait $S$. Ceci contredirait l'hypoth\`ese $\mathcal{O}(C_s)\notin\Pic(X)$. L'ensemble $\Sigma'$ est donc dense dans $S$.

Si $s\in\Sigma'$, la th\'eorie des d\'eformations montre que le lieu sur $U$ o\`u le fibr\'e en droites $\mathcal{O}(C_s)$ se d\'eforme est un diviseur irr\'eductible : on le note $B_s$ et on note $\bar{\eta}_s$ un point g\'en\'erique g\'eom\'etrique de $B_s$. Par la proposition \ref{defo}\footnote{L'hypoth\`ese de rigidit\'e est v\'erifi\'ee car une surface K3 en caract\'eristique nulle n'est pas unir\'egl\'ee.}, $f_s$ se d\'eforme le long de $B_s$ : il existe une courbe rationnelle $\widetilde{C}_s\subset\mathcal{X}_{\bar{\eta}_s}$ dont $C_s$ est une sp\'ecialisation. Les $(B_s)_{s\in\Sigma'}$ forment une collection infinie de diviseurs : dans le cas contraire, par densit\'e de $\Sigma'$ dans $S$, un des $B_s$ contiendrait $S$, et cela contredirait l'hypoth\`ese que $\mathcal{O}(C_s)\notin\Pic(X)$. Par cons\'equent,  la r\'eunion des $(B_s)_{s\in\Sigma'}$ est dense dans $U$. On peut alors appliquer le th\'eor\`eme \ref{nodale} : il existe $s\in\Sigma'$ tel que le syst\`eme lin\'eaire $|\mathcal{A}_{\bar{\eta}_s}|$ contienne une courbe rationnelle nodale $R$, qui est bien s\^ur ample. De plus, par le th\'eor\`eme \ref{unecourbe}, si $m\gg 0$, il existe un diviseur $D\in|\mathcal{A}_{\bar{\eta}_s}^{\otimes m}(-\widetilde{C}_s)|$ dont toutes les composantes irr\'eductibles sont des courbes rationnelles.

On peut alors appliquer la proposition \ref{rigide} : il existe une courbe stable $[g]\in\overline{\mathcal{M}}_0(\mathcal{X}_{\bar{\eta}_s})$ dont le cycle associ\'e est $\widetilde{C}_s+D+kR$, et qui est g\'en\'eriquement non ramifi\'ee et rigide.

 Par la proposition \ref{defo}, on peut choisir une composante irr\'eductible $M$ de $\overline{\mathcal M}_0(\mathcal{X}/U)$ contenant $[g]$ qui domine $U$. Celle-ci est de plus propre sur $U$, par propret\'e des champs de modules d'applications stables. D'une part, le choix de $\widetilde{C}_s$ et la propret\'e de $M\to U$ montre qu'il est possible de sp\'ecialiser $g$ en une courbe stable $g_s$ sur $\mathcal{X}_s$ dont le cycle associ\'e contient $C_s$. D'autre part, on peut choisir une courbe stable $g_{\bar{\xi}}$ sur $\mathcal{X}_{\bar{\xi}}$ qui est une g\'en\'erisation de $g$ dans $M$. On peut alors choisir une sp\'ecialisation $g_{\bar{\eta}}$ de $g_{\bar{\xi}}$ sur~$\mathcal{X}_{\bar{\eta}}$ dont $g_s$ est une sp\'ecialisation. Le cycle associ\'e \`a $g_{\bar{\eta}}$ contient alors une composante irr\'eductible de degr\'e $>N$, car tel est le cas pour sa sp\'ecialisation $g_s$.

On a donc construit une courbe rationnelle de degr\'e $>N$ sur $X$.
Comme $N$ \'etait arbitraire, on a montr\'e l'existence d'une infinit\'e de courbes rationnelles dans $X$.
\qed

\subsection{Surfaces K3 de nombre de Picard impair}\label{parr4}

Nous pouvons \`a pr\'esent expliquer la preuve du th\'eor\`eme \ref{BHTLL2} quand $\rho(X)$ est impair.

\begin{theo}\label{th2}
Soit $X$ une surface K3 sur un corps alg\'ebriquement clos de carac\-t\'eristique $\neq 2$ telle que $\rho(X)$ est impair. Alors $X$~contient
une infinit\'e de courbes rationnelles.
\end{theo}

\noindent{\sc Preuve} --- 
La preuve du th\'eor\`eme \ref{th1} ci-dessus s'adapte : on utilise toujours des arguments de sp\'ecialisation et d\'eformation, mais en caract\'eristique mixte. Contentons-nous d'indiquer les quelques modifications \`a apporter.

  On consid\`ere \`a pr\'esent une carte \'etale $U$ du champ de modules $\mathfrak{M}_{2d,\mathbb{Z}[\frac{1}{2}]}$ en carac\-t\'eristique mixte, et on conserve les autres notations. En tous les points ferm\'es $s\in S$, le nombre de Picard g\'eom\'etrique cro\^it  : cela r\'esulte de l'hypoth\`ese sur $\rho(X)$ et de la conjecture de Tate pour les surfaces K3 sur les corps finis (via le corollaire \ref{cororho} (iii)). Cela remplace l'argument de th\'eorie de Hodge utilis\'e dans la preuve du th\'eor\`eme \ref{th1}.

  On d\'efinit alors $\Sigma$ comme l'ensemble des points ferm\'es non supersinguliers de $S$. Si $X$ est de caract\'eristique finie, la conjecture d'Artin (corollaire \ref{cororho} (ii)) montre que $X$ ne peut pas \^etre supersinguli\`ere, et $\Sigma$ est dense dans $S$ par fermeture du lieu supersingulier. Si $X$ est de caract\'eristique nulle, un th\'eor\`eme\footnote{Ce th\'eor\`eme permet m\^eme de ne conserver dans $\Sigma$ que des surfaces K3 ordinaires. Dans ce cas, les travaux de Nygaard \cite{Nygaard} sur la conjecture de Tate sont donc suffisants.} de Joshi et Rajan \cite[Theorem 6.6.2]{JR} et Bogomolov et Zarhin \cite{BZ} assure que $\Sigma$ est dense dans $S$.

  Enfin, il faut utiliser une variante en caract\'eristique mixte du th\'eor\`eme \ref{defo} (voir la remarque \ref{deforem} et \cite[Theorem 18]{BHT}) pour d\'eformer $f_s$, qui n\'ecessite deux pr\'ecautions. D'une part, pour v\'erifier l'hypoth\`ese de rigidit\'e dans cet \'enonc\'e, nous avons utilis\'e le fait que les surfaces K3 en caract\'eristique nulle ne sont pas unir\'egl\'ees. En caract\'eristique finie, ce n'est pas vrai en g\'en\'eral, mais c'est le cas pour les surfaces K3 non supersinguli\`eres \cite[Corollary 2]{Shiodaex}\footnote{Il est montr\'e ici qu'une surface K3 unirationnelle est supersinguli\`ere, mais la m\^eme preuve fonctionne si elle est unir\'egl\'ee.} : c'est la raison pour laquelle nous avons modifi\'e la d\'efinition de $\Sigma$. D'autre part, le lieu $D_s$ dans $U$ o\`u $\mathcal{O}(C_s)$ se d\'eforme est toujours un diviseur, est plat sur $\mathbb{Z}$, mais n'est pas n\'ecessairement irr\'eductible si le fibr\'e en droites $\mathcal{O}(C_s)$ n'est pas primitif. Le dernier argument de la preuve du th\'eor\`eme \ref{defo} ne s'\'etend pas tel quel \`a cette situation : il permet seulement de montrer que $f_s$ se d\'eforme le long d'une des composantes irr\'eductibles de $D_s$. On d\'efinit alors $B_s$ comme \'etant l'une de ces composantes irr\'eductibles.
  
  Le reste de la preuve fonctionne sans modifications.
  \qed

\subsection{Limites de la strat\'egie}\label{parr5}

Des techniques diff\'erentes ont permis de montrer la conjecture \ref{conjcourbes} pour les surfaces K3 elliptiques en caract\'eristique nulle \cite{BT,Hassettsurvey}. Comme une surface K3 de nombre de Picard $\geq 5$ est automatiquement elliptique, le th\'eor\`eme \ref{BHTLL2} implique :

\begin{coro}
Soit $X$ une surface K3 sur $\ovQ$ telle que $\rho(X)\notin\{2,4\}$. Alors $X$contient
une infinit\'e de courbes rationnelles.
\end{coro}

Par cons\'equent, les surfaces K3 en caract\'eristique nulle pour lesquelles la conjecture \ref{conjcourbes} n'est pas connue sont d\'efinies sur $\ovQ$ et ont nombre de Picard $2$ ou $4$. Si l'on souhaitait utiliser la strat\'egie de la preuve du th\'eor\`eme \ref{BHTLL2} pour une telle surface K3 $X$ d\'efinie sur un corps de nombres $K$, il faudrait conna\^itre l'existence d'une infinit\'e de places $\mathfrak{p}$ de $K$ telles que le nombre de Picard g\'eom\'etrique de $X_{\mathfrak{p}}$ soit sup\'erieur \`a celui de $X$ (et telles que $X_{\mathfrak{p}}$ ne soit pas supersinguli\`ere). 

C'est une question ouverte. Sa difficult\'e tient au fait, remarqu\'e par Charles \cite[Theorem 1]{CharlesPic}, que l'ensemble de ces places peut \^etre de densit\'e nulle. L'analyse de Charles d\'ecrit compl\`etement ce ph\'enom\`ene. Elle permet de donner des exemples de surfaces K3 sur $\ovQ$ de nombre de Picard $2$ ou $4$ pour lesquelles la strat\'egie s'applique, et qui contiennent donc une infinit\'e de courbes rationnelles \cite[Corollary 4]{CharlesPic}.

\bigskip

Enfin, la strat\'egie de Bogomolov, Hassett et Tschinkel repose de mani\`ere essentielle sur l'existence de sp\'ecialisations, qui permettent de s'assurer que l'on construit bien une infinit\'e de courbes rationnelles distinctes. En particulier, cette strat\'egie ne permet de rien dire sur les surfaces K3 d\'efinies sur $\ovFp$. En un sens, ce cas est le plus important. 

De plus, cette strat\'egie ne donne pas de contr\^ole sur les courbes rationnelles construites. Par exemple, elle ne permet pas de construire une courbe rationnelle \'evitant un point fix\'e ou appartenant \`a un multiple d'un syst\`eme lin\'eaire fix\'e. Elle laisse \'egalement compl\`etement ouverte la question suivante de Bogomolov : 
\begin{ques}
Soit $X$ une surface K3 sur $\ovQ$ ou $\ovFp$. Passe-t-il une courbe rationnelle par tout point ferm\'e de $X$ ? 
\end{ques}
Sur $\ovQ$, cette question est qualifi\'ee pudiquement de \og logical possibility\fg \ dans \cite{BT2}, et on n'en conna\^it ni exemple ni contre-exemple. Sur $\ovFp$, Bogomolov et Tschinkel \cite[Theorem 1.1]{BT2}  ont obtenu le cas particulier suivant, montrant que l'\'enonc\'e est plausible :
\begin{theo}
Soit $X$ une surface de Kummer sur $\ovFp$. Alors il passe une courbe rationnelle par tout point ferm\'e de $X$.
\end{theo}


\end{document}